\documentclass[12pt,oneside,french]{article}
\usepackage{amssymb,amsmath,babel}
\newtheorem{Theoreme}{Th\'eor\`eme}
\newtheorem{Lemme}{Lemme}
\newtheorem{Proposition}{Proposition}

\newcommand{\NN}{\mathbb N}
\newcommand{\PP}{\mathbb P}
\newcommand{\ZZ}{\mathbb Z}
\newcommand{\QQ}{\mathbb Q}
\newcommand{\RR}{\mathbb R}
\newcommand{\CC}{\mathbb C}

\newcommand{\sqm}[4]{\displaystyle{
\left({#1 \atop #3}{#2 \atop #4}\right)}}

\def\und{\underline}
\def\SL{{\bf SL}}

\begin{document}

\begin{center}

\LARGE{Lemmes de multiplicit\'e associ\'es aux \\ groupes triangulaires de Riemann-Schwarz.}

\vspace{15pt}

\Large{Federico Pellarin}

\end{center}

\vspace{20pt}

\noindent {\footnotesize {\bf English abstract.} In this article, we extend a multiplicity estimate of Nesterenko, valid for 
quasi-modular forms associated to $\SL_2(\ZZ)$, to non-holomorphic
quasi-modular forms associated to certain co-compact Riemann-Schwarz triangular subgroups of $\SL_2(\CC)$.}

\section{Introduction.\label{section:introduction}}
Soient $E_2(z),E_4(z),E_6(z)$
les d\'eveloppements de Fourier complexes des s\'eries d'Eisenstein classiques de poids $2,4,6$, convergents pour $|z|<1$.

Nesterenko a demontr\'e que pour tout nombre complexe $q$ tel que $0<|q|<1$, le corps $\QQ(q,E_2(q),E_4(q),E_6(q))$
a degr\'e de transcendance au moins $3$ (voir \cite{Nesterenko:Article},
\cite{Nesterenko:Introduction3} et \cite{Nesterenko:Introduction10}). L'ingredient cl\'ef
de sa preuve est le lemme de multiplicit\'e ci-dessous (cf. th\'eor\`eme 2.3 p. 33 de
\cite{Nesterenko:Introduction3}).

\begin{Theoreme}[Nesterenko] Il existe une constante $c_1>0$ avec la propri\'et\'e suivante.
Soit $P$ un polyn\^{o}me non nul de $\CC[X_1,X_2,X_3,X_4]$, de degr\'e total au plus $N$.
Alors, la fonction $F(z)=P(z,E_2(z),E_4(z),E_6(z))$ s'annule en $z=0$ avec une multiplicit\'e
au plus $c_1N^4$.
\label{theoreme:estimate_nesterenko}\end{Theoreme}
Le but de ce texte est de donner des g\'en\'eralisations du th\'eor\`eme \ref{theoreme:estimate_nesterenko}. 
Soient
$\alpha,\beta,\gamma$ des inverses d'entiers naturels
non nuls tels que~:
\begin{equation}
\gamma>\alpha+\beta,\quad
1>\gamma>\beta>\alpha>0.\label{eq:conditionsabc}\end{equation}
On consid\`ere l'\'equation diff\'erentielle hyperg\'eom\'etrique complexe~:
\begin{equation}
z(1-z)\frac{d^2V}{d z^2}+(\gamma-(\alpha+\beta+1)z)\frac{dV}{dz}-\alpha\beta V=0.\label{eq:hypergeometrique}\end{equation}
 Le groupe de monodromie projective de l'\'equation 
(\ref{eq:hypergeometrique}) s'identifie \`a un sous-groupe $\Omega$ discret infini de $\SL_2(\CC)$ (un sous-groupe Fuchsien de
premi\`ere esp\`ece,
cf. \cite{Ford:Automorphic}, chapitre 3), qui agit  discontinuement sur le disque $B=\{z\in\CC\mbox{ tel que
}|z|<1\}$ par transformations homographiques~:
\[\phi=\sqm{a}{b}{c}{d}\in\Omega,\quad\phi(\xi):=\frac{a\xi+b}{c\xi+d}.\]
Il existe un domaine fondamental
${\cal T}$ pour l'action de $\Omega$ sur $B$, dont l'adh\'erence topologique  est un triangle hyperbolique compact, de
sommets $s_0,s_1,s_\infty$, et dont les angles aux sommets sont \'egaux
\`a $\pi(1-\gamma),\pi(\gamma-\alpha-\beta),\pi(\beta-\alpha)$ (cf. \cite{Yoshida:Love} chapitre 3 ou \cite{Yoshida:Fuchsian}
chapitre 5). Posons 
\[B^*=\bigcup_{\phi\in\Omega}\phi({\cal T}-\{s_0,s_1,s_\infty\}).\]
Ainsi, $B^*$ est le disque $B$ priv\'e de l'ensemble $E$ dont les \'el\'ements sont 
tous les points $\xi$ tels qu'il existe
$\gamma\in\Omega$ avec \[\gamma(\xi)\in\{s_0,s_1,s_\infty\}.\]
On peut construire (cf. \cite{Zudilin:Hypergeometric}) trois fonctions $Y_0(t),Y_1(t),Y_2(t)$ holomorphes
dans $B^*$ et alg\'ebri\-quement ind\'ependantes, qui satisfont des 
relations d'automorphie de poids $2$ par rapport \`a l'action de $\Omega$ sur $B$. 
Plus pr\'eci\-sement, si $\phi=\sqm{a}{b}{c}{d}\in\Omega$, alors on a~:
\begin{equation}
Y_i(\phi(t))=(ct+d)^2Y_i(t)+\frac{1}{\pi{\rm i}}c(ct+d),\quad i=0,1,2.\label{eq:automorphes}\end{equation}
Les conditions ci-dessus d\'eterminent $Y_0(t),Y_1(t),Y_2(t)$ \`a une constante multiplicative pr\`es,
et \`a permutation pr\`es.

De plus, les conditions (\ref{eq:conditionsabc}) d\'eterminent de mani\`ere unique 
un entier positif non nul $q$ tel que, au voisinage de $s_0,s_1,s_\infty$, on ait les d\'eveloppements
en s\'erie de Puiseux-Laurent convergents (dans la suite appel\'es plus simplement d\'evelop\-pements en s\'erie de Puiseux)~:
\begin{equation}
Y_i(t) = \sum_{k=-q}^\infty v_{i,k}(t-s_j)^{k/q},\quad\mbox{avec }i=0,1,2\mbox{ et }j=0,1,\infty,
\label{eq:Puiseux_laurent}\end{equation}
o\`u $v_{i,k}\in\CC$ et $v_{i,-q}\not=0$ pour tout $i=0,1,2$.

Les cinq fonctions $t,e^t,Y_0(t),Y_1(t),Y_2(t)$ sont alg\'ebriquement ind\'epen\-dantes (cf. 
\cite{Nishioka:Conjecture}, voir aussi \cite{Nesterenko:Differential}, ou \cite{Zudilin:Hypergeometric},
proposition 5). 

Tout polyn\^{o}me $F$ en $t,e^t$ et les fonctions $Y_i(t)$, ($i=0,1,2$) \`a coefficients 
complexes admet un d\'eveloppement en s\'erie
de Puiseux convergent au voisinage de tout \'el\'ement $\xi\in{\cal T}$:
\[F(t)=\sum_{k=s}^\infty v_k(t-\xi)^{k/q},\] avec $s\in\ZZ$ et $v_{s}\not=0$ (en particulier, $F(\xi+t^q)$
est une fonction m\'eromorphe au voisinage de $\xi$).
On d\'efinit alors la multiplicit\'e de $F$ en $\xi$ par:
\[\mbox{ord}_{\xi}(F)=\frac{s}{q}.\]
Si $\xi\not=s_0,s_1,s_\infty$, alors $\mbox{ord}_{\xi}(F)\in\NN$, autrement, 
$\mbox{ord}_{\xi}(F)\in\ZZ/q$ peut \^{e}tre n\'egatif.

Le but de ce texte est de
d\'emontrer des estimations de multiplicit\'e g\'en\'e\-ralisant le th\'eor\`eme \ref{theoreme:estimate_nesterenko}
aux \'el\'ements de \[\CC[t,e^t,Y_0(t),Y_1(t),Y_2(t)],\] que nous r\'eunissons
dans le th\'eor\`eme suivant.

\begin{Theoreme}
Il existe une constante $c_2>0$, d\'ependant seulement de $\alpha,\beta,\gamma$, avec les
propri\'et\'es suivantes. Soit $\xi\in B$, soit $P$ un polyn\^{o}me non nul de
$\CC[X_1,X_2,X_3,X_4,X_5]$, posons~: 
\begin{eqnarray*}
M_1& = & \min\{\deg_{X_1}(P),\deg_{X_2}(P)\}+1\\
M_2& = &\max\{\deg_{X_1}(P),\deg_{X_2}(P)\}+\max\{\deg_{X_3}(P),\deg_{X_4}(P),\deg_{X_5}(P)\}.
\end{eqnarray*} Alors la fonction
$F(z)=P(t,e^t,Y_0(t),Y_1(t),Y_2(t))$ satisfait
\begin{equation}
\mbox{ ord}_{\xi}(F)\leq c_2M_1M_2^4.\label{eq:est1}\end{equation}

En particulier, si $P$ est de degr\'e total au plus $N$, alors 
la fonction $F(z)$ satisfait
\begin{equation}
\mbox{ ord}_{\xi}(F)\leq c_2(2N)^5.\label{eq:est2}\end{equation} 

Si de plus $P$ ne d\'epend pas de la variable $X_1$ (ou si $P$ ne d\'epend pas de la variable $X_2$),
alors
\begin{equation}\mbox{ ord}_{\xi}(F)\leq c_2(2N)^4.\label{eq:est3}\end{equation}
\label{theorem:estimate2}\end{Theoreme}
Ce th\'eor\`eme g\'en\'eralise le th\'eor\`eme \ref{theoreme:estimate_nesterenko} de Nesterenko. Il contient aussi
un lemme de multiplicit\'e de Bertrand~: lemme 3 p. 348 de \cite{Bertrand:Theta}. 

Soit ${\cal A}$ un anneau muni d'une d\'erivation
$\delta:{\cal A}\rightarrow{\cal A}$~; nous dirons que le couple $({\cal A},\delta)$ est un 
{\em anneau diff\'erentiel}.

La d\'emonstration du th\'eor\`eme \ref{theoreme:estimate_nesterenko} utilise en profondeur le fait que l'anneau
\[{\cal Y}_1=\CC[z,E_2(z),E_4(z),E_6(z)],\] muni de la d\'erivation $z(d/dz)$, est un anneau diff\'erentiel 
(cf. th\'eor\`eme 5.3 de \cite{Lang:Modular}).

De mani\`ere analogue, l'anneau \[{\cal A}=\CC[t,e^t,Y_0(t),Y_1(t),Y_2(t)],\] muni de la d\'erivation
$\delta=d/dt$, 
poss\`ede une structure d'anneau diff\'erentiel (cf. lemme \ref{lemme:relations_explicites} de ce texte). 

Rappelons ici qu'un id\'eal ${\cal P}$ d'un anneau diff\'erentiel $({\cal A},\delta)$ est dit
{\em $\delta$-stable} si pour tout $x\in{\cal P}$ on a $\delta x\in{\cal P}$.
Nesterenko d\'emontre le r\'esultat qui suit, d'inter\^{e}t ind\'ependant (cf. proposition 5.1 p. 161 de
\cite{Nesterenko:Introduction3}), et indispensable dans la preuve du th\'eor\`eme \ref{theoreme:estimate_nesterenko}.
\begin{Proposition}
Soit ${\cal P}$ un id\'eal premier non nul et $z(d/dz)$-stable de ${\cal Y}_1$,
tel que pour tout $F\in{\cal P}$ on ait $F(0)=0$. Alors $z(E_4(z)^3-E_6(z)^2)\in{\cal P}$.
\label{proposition:classif_nest}\end{Proposition}

Pour d\'emontrer le th\'eor\`eme \ref{theorem:estimate2}, nous utiliserons la propri\'et\'e diff\'e\-ren\-tielle de 
$({\cal A},\delta)$ d\'ecrite par la proposition suivante.
\begin{Proposition} Il existe un \'el\'ement non nul $\kappa\in{\cal A}$
tel que $\kappa\in{\cal P}$ pour tout id\'eal premier non nul $\delta$-stable ${\cal P}$ de ${\cal A}$.
\label{proposition:prop_ramanujan1}\end{Proposition}
La d\'emonstration de la proposition \ref{proposition:classif_nest} donn\'ee par Nesterenko 
utilise de ma\-ni\`ere essentielle la non compacit\'e du quotient $\SL_2(\ZZ)\backslash{\cal H}$,
o\`u ${\cal H}$ est le demi-plan sup\'erieur (existence de d\'eveloppements en s\'erie de Fourier de
formes modulaires et quasi-modulaires au voisinage de {\em pointes g\'eom\'etriques}).

Pour d\'emontrer la proposition \ref{proposition:prop_ramanujan1}
nous suivons de pr\`es nombreuses techniques d\'ej\`a introduites par Nesterenko dans \cite{Nesterenko:Introduction3},
\cite{Nesterenko:Introduction10} (voir le plan de l'article ci-dessous). 
Cependant, dans un point au moins, les m\'ethodes de
Nesterenko ne s'\'etendent pas aux fonctions $Y_0,Y_1,Y_2$~: sous les hypoth\`eses du
th\'eor\`eme \ref{theorem:estimate2},
$\Omega\backslash B$ est compact, contrairement \`a $\SL_2(\ZZ)\backslash{\cal H}$.

\subsection{Plan de l'article et structure des d\'emonstrations.}

La d\'emonstration de la proposition \ref{proposition:prop_ramanujan1} occupe
tout le paragraphe \ref{section:systemes}. En absence de pointes
dans le quotient $\Omega\backslash B$, nous avons effectu\'e une param\'etri\-sation de $Y_0,Y_1,Y_2$ qui nous a
induit \`a consid\'erer des fonctions hyperg\'eo\-m\'etriques. Dans ce cadre, nous faisons jouer 
aux points $0,1,\infty\in\PP_1(\CC)$
un r\^{o}le similaire \`a celui de la pointe \`a l'infini de $\SL_2(\ZZ)\backslash{\cal H}$, ce qui nous permet
de g\'en\'eraliser le lemme 5.2 p. 161 de \cite{Nesterenko:Introduction10} (proposition 
\ref{lemme:ideaux_principaux_hypergeo}
de ce texte). De
cette mani\`ere nous d\'emontrons la proposition \ref{proposition:prop_ramanujan1} dans le cas particulier
o\`u ${\cal P}$ est 
un id\'eal principal.

Pour traiter le cas des id\'eaux non principaux dans la proposition \ref{proposition:prop_ramanujan1},
nous avons introduit une technique compl\`etement alg\'ebrique qui utilise 
les {\em crochets de Rankin} issus de la th\'eorie des formes quasi-modulaires. 
Cette technique permet aussi de raffiner et g\'en\'eraliser les r\'esultats
du paragraphe 5 pp. 162-165 de \cite{Nesterenko:Introduction10}~: par exemple, on peut facilement supprimer 
l'hypoth\`ese $F(0)=0$ dans la proposition \ref{proposition:classif_nest}~: voir \cite{Pellarin:Quasimodular}.

Passons aux r\'esultats du paragraphe \ref{section:demo_th}.
Nous y supposerons que $({\cal A},\delta)$ soit plus g\'en\'eralement un 
anneau diff\'erentiel de fonctions $F$ admettant en tout point $\xi$ d'un certain domaine
de $\CC$,
un d\'eveloppement en s\'erie de Puiseux convergent~:
\[F(t)=\sum_{k=-s}^\infty v_k(t-\xi)^{k/q},\] pour un certain entier $q>0$ ne d\'ependant pas de $\xi$.

Cet anneau diff\'erentiel $({\cal A},\delta)$ satisfait aussi une \og propri\'et\'e $D$\fg 
(cette propri\'et\'e est tr\`es proche d'une propri\'et\'e homonyme dans 
\cite{Nesterenko:Introduction10}~: d\'efinition 1.2 p. 150).
Sous toutes ces hypoth\`eses, on obtient un lemme de multiplicit\'e~: la
proposition 
\ref{proposition:general}.

Si la propri\'et\'e $D$ est v\'erifi\'ee pour l'anneau diff\'erentiel $({\cal A},\delta)$
de la proposition \ref{proposition:prop_ramanujan1}, alors on obtient le th\'eor\`eme
\ref{theorem:estimate2} comme corollaire  de la proposition \ref{proposition:general}.
Or, nous v\'erifierons que la propri\'et\'e d\'ecrite par la proposition \ref{proposition:prop_ramanujan1}
(que nous appelons \og propri\'et\'e de Ramanujan\fg ) implique la propri\'et\'e $D$ pour $({\cal A},\delta)$.

Les techniques de d\'emonstration de la proposition \ref{proposition:general}
sont essentiellement les m\^{e}mes que celles de \cite{Nesterenko:Introduction10}, bien que nous trouvions 
pr\'eferable de suivre l'approche
de \cite{Bosser:Independance}. Dans le paragraphe \ref{section:demo_th} nous d\'etaillerons autant que possible des
d\'emonstrations pourtant tr\`es voisines de celles des deux r\'eferences ci-dessus, uniquement dans le souci de traiter un cas
apparemment plus g\'en\'eral. En effet, les fonctions $Y_0,Y_1,Y_2$ ne sont pas holomorphes dans ${\cal T}$, \`a cause des
singularit\'es  en $s_0,s_1,s_\infty$. Or dans
\cite{Bosser:Independance} et \cite{Nesterenko:Introduction10}, on travaille uniquement avec des fonctions 
holomorphes.

\section{Anneaux diff\'erentiels.\label{section:systemes}}

Dans ce paragraphe nous d\'emontrons la proposition \ref{proposition:prop_ramanujan1}.
Le polyn\^{o}me $\kappa$ sera explicitement d\'etermin\'e
(cf. proposition \ref{proposition:stables} du paragraphe \ref{section:ideaux_stables}), m\^{e}me si cette donn\'ee suppl\'ementaire
n'est pas utilis\'ee dans la d\'emonstration du th\'eor\`eme \ref{theorem:estimate2}. Nous allons commencer par une \'etude
diff\'erentielle des fonctions hyperg\'eom\'etriques.

Consid\'erons les nombres rationnels~:
\begin{eqnarray*}
a & = & \gamma(1-\alpha-\beta)+2\alpha\beta\\
b & = & (\alpha+\beta)(\gamma-\alpha-\beta)+2\alpha\beta-\gamma+1\\
c & = & \gamma(\alpha+\beta-\gamma+1)-2\alpha\beta.
\end{eqnarray*}
Un calcul direct permet de v\'erifier que 
les fonctions 
\begin{eqnarray*}
u_0(z) & = & z^{\gamma/2}(1-z)^{(\alpha+\beta-\gamma+1)/2}{}_2F_1(a,b;c;z),\\
u_1(z) & = & z^{1-\gamma/2}(1-z)^{(\alpha+\beta-\gamma+1)/2}{}_2F_1(\alpha-\gamma+1,\beta-\gamma+1;2-\gamma;z),
\end{eqnarray*}
analytiques dans $\CC-(\RR_{\leq 0}\cup\RR_{\geq 1})$, sont $\CC$-lin\'eairement ind\'ependantes et 
sa\-tisfont l'\'equation diff\'erentielle (cf. \cite{Ford:Automorphic}, \'equation (21) p. 290)
\begin{equation}
\frac{d^2U}{dz^2}+\left(\frac{a}{4z^2}+\frac{b}{4(z-1)^2}+\frac{c}{4z^2(z-1)^2}\right)U=0
\label{eq:nouvelle}.\end{equation}
Pour $z\in\CC-(\RR_{\leq 0}\cup\RR_{\geq 1})$, posons~:
\begin{eqnarray*}
y_0  & = & u_0 u_0' ,\\
y_1  & = & u_0 u_0' -\frac{u_0^2}{z},\\
y_2  & = & u_0 u_0' -\frac{u_0^2}{z-1},\\
\tau  & = & \frac{u_1 }{u_0 },\\
q  & = & e^{\tau }.
\end{eqnarray*}
La fonction $\tau(z)$ est localement inversible sur son domaine d'holomorphie, car le
wronskien 
$\det\sqm{u_1}{u_0}{u_1'}{u_0'}$ de $u_1,u_0$ y est constant, \'egal \`a $1-\gamma$, donc non nul. Dans la suite, nous posons
$w=1-\gamma$.

Soit ${\cal H}$ le demi-plan sup\'erieur complexe. On peut montrer que la fonction $\tau$ d\'efinit un isomorphisme analytique
\[\tau:{\cal H}\rightarrow{\cal T}^\circ,\] o\`u ${\cal T}^\circ$ d\'esigne le triangle ${\cal T}$
priv\'e de son bord topologique; on v\'erifie aussi $\tau(i)=s_i$ avec $i=0,1,\infty$.

Soit $\zeta:{\cal T}^\circ\rightarrow{\cal H}$ la fonction analytique r\'eciproque de $\tau$~:
elle admet un prolongement analytique \`a $B^*$ et d\'efinit sur $B$ une fonction alg\'ebrique sur le corps
des fonctions modulaires Fuchsiennes associ\'ees \`a $\Omega$, c'est-\`a-dire les fonctions $f:B\rightarrow\PP_1(\CC)$
satisfaisant~:
\[f(\phi(t))=f(t),\quad\phi\in\Omega\]
(cf. \cite{Ford:Automorphic}, chapitre 10). Posons
\[Y_i(t)=y_i(\zeta(t))\] pour
$i=0,1,2$. On v\'erifie que~:
\begin{eqnarray*}
Y_0(t)& = & w\frac{d}{dt}\log (u_0),\\
Y_1(t)& = & w\frac{d}{dt}\log \left(\frac{u_0}{z}\right),\\
Y_2(t)& = & w\frac{d}{dt}\log \left(\frac{u_0}{z-1}\right),\\
\end{eqnarray*}
et c'est facile de montrer, \`a partir de ces formules, que ces fonctions satisfont les relations d'automorphie
(\ref{eq:automorphes}) (cf. \cite{Zudilin:Hypergeometric}, proposition 3). Ce sont les
fonctions $Y_0,Y_1,Y_2$ du paragraphe \ref{section:introduction}. 

Soit $D$ la d\'erivation d\'efinie sur l'anneau des fonctions holomorphes sur
$\CC-(\RR_{\leq 0}\cup\RR_{\geq 1})$ par \[DX=\displaystyle{u_0^2\frac{dX}{dz}}.\] L'anneau
\[{\cal N}=\CC[\tau,q,y_0,y_1,y_2]\] est stable pour la d\'erivation $D$, comme le montre le lemme qui suit.
\begin{Lemme} Soit $L=(1/4)(a(y_0-y_1)^2+b(y_0-y_2)^2+c(y_1-y_2)^2)$. Les relations suivantes sont
satisfaites. \[D\tau=w,\quad Dq=wq,\quad Dy_i=y_i^2-L,\quad \mbox{ pour }i=0,1,2.\]
\label{lemme:relations_explicites}\end{Lemme}
\noindent {\bf D\'emonstration}. On a 
\[D\tau = u_0^2\tau'= u_0^2\left(\frac{u_1'u_0-u_1u_0'}{u_0^2}\right)
= W(u_1,u_0)= w.\]
De m\^{e}me, $Dq=u_0^2\tau'q=wq$.

Calculons maintenant $Dy_0,Dy_1,Dy_2$. On v\'erifie facilement l'identit\'e~:
\begin{eqnarray*}
L & = & \frac{1}{4}u_0^4\left(\frac{a}{4z^2}+\frac{b}{4(z-1)^2}+\frac{c}{4z^2(z-1)^2}\right).
\end{eqnarray*}
On a, en utilisant (\ref{eq:nouvelle})~:
\begin{eqnarray*}
Dy_0 & = & u_0^2y_0'\\
& = & u_0^2(u_0u_0')'\\
& = & u_0^2(u_0'{}^2+u_0u_0'')\\
& = & (u_0u_0')^2+u_0^3u_0''\\
& = & y_0^2-u_0^4\left(\frac{a}{4z^2}+\frac{b}{4(z-1)^2}+\frac{c}{4z^2(z-1)^2}\right)\\
& = & y_0^2-L.
\end{eqnarray*}
De m\^{e}me,
\begin{eqnarray*}
Dy_1 & = & u_0^2y_1'\\
& = & u_0^2(u_0u_0'-u_0^2z^{-1})'\\
& = & u_0^2(u_0'{}^2+u_0u_0''-2u_0u_0'z^{-1}+u_0^2z^{-2})\\
& = & (u_0u_0')^2-2u_0^3u_0'z^{-1}+u_0^4z^{-2}+u_0^3u_0''\\
& = & (u_0u_0'-u_0^2z^{-1})^2+u_0^3u_0''\\
& = & y_1^2-L,
\end{eqnarray*}
et pour terminer,
\begin{eqnarray*}
Dy_2 & = & u_0^2y_2'\\
& = & u_0^2(u_0u_0'-u_0^2(z-1)^{-1})'\\
& = & u_0^2(u_0'{}^2+u_0u_0''-2u_0u_0'(z-1)^{-1}+u_0^2(z-1)^{-2})\\
& = & (u_0u_0')^2-2u_0^3u_0'(z-1)^{-1}+u_0^4(z-1)^{-2}+u_0^3u_0''\\
& = & (u_0u_0'-u_0^2(z-1)^{-1})^2+u_0^3u_0''\\
& = & y_2^2-L.
\end{eqnarray*} Donc la d\'erivation $D$ laisse stable l'anneau ${\cal N}$. On d\'eduit du lemme
\ref{lemme:relations_explicites} que
l'isomorphisme local $\zeta$ induit un isomorphisme d'anneaux diff\'erentiels $({\cal N},D)\cong({\cal
A},\delta)$, o\`u $\delta=w(d/dt)$ et ${\cal A}=\CC[t,e^t,Y_0,Y_1,Y_2]$.
\subsection{Propri\'et\'es asymptotiques en $0,1,\infty$.}
Nous \'etudions les s\'eries de Puiseux des fonctions $u_0^2,y_0,y_1,y_2$ au voisinage de
$0,1,\infty\in\PP_1(\CC)$. Tout d'abord, introduisons quelques notations. Soit $p\in\NN$ le plus petit
d\'enominateur commun des nombres rationnels $\alpha,\beta,\gamma$. Nous allons travailler
dans les anneaux de s\'eries de Puiseux
formelles~:
\begin{eqnarray*}
{\cal T}_0 & = & \CC[[z^{1/p}]],\\
{\cal T}_1 & = & \CC[[(1-z)^{1/p}]],\\
{\cal T}_\infty & = & \CC[[(-z)^{-1/p}]].
\end{eqnarray*}
Nous notons ${\cal M}_0,{\cal M}_1,{\cal M}_\infty$ les id\'eaux maximaux des 
s\'eries de ${\cal T}_0,{\cal T}_1,{\cal T}_\infty$ sans terme constant, et 
${\cal Q}_0,{\cal Q}_1,{\cal Q}_\infty$ les corps des fractions de ${\cal T}_0,{\cal T}_1,{\cal T}_\infty$.

Une s\'erie $\Sigma_0$ de ${\cal T}_0$ est dite {\em localement convergente}
si elle s'\'ecrit~:
\[\Sigma_0(z)=\sum_{i=0}^{p-1}f_i(z)z^{i/p},\] avec $f_i\in\CC[[z]]$ convergente au voisinage de $0$.
De m\^{e}me, nous dirons qu'une s\'erie $\Sigma_1$ de ${\cal T}_1$ est {\em localement convergente}
si elle s'\'ecrit~:
\[\Sigma_1(z)=\sum_{i=0}^{p-1}f_i(z)(1-z)^{i/p},\] avec $f_i\in\CC[[1-z]]$ convergente au voisinage de $0$.
De mani\`ere \'equivalente, on peut observer que si $\Sigma_1\in{\cal T}_1$, alors $\Sigma_1(1-z)\in{\cal T}_0$~;
donc $\Sigma_1$ est localement convergente si $\Sigma_1(1-z)$ l'est.

Soit $\Sigma_\infty$ une s\'erie de ${\cal T}_\infty$. Alors $\Sigma_\infty(-z^{-1})\in {\cal T}_0$.
Nous dirons que $\Sigma_\infty$ est {\em localement convergente} si $\Sigma_\infty(-z^{-1})$ est localement 
convergente.

Ces d\'efinitions permettent de caract\'eriser \'egalement les s\'eries localement convergentes de 
${\cal Q}_0,{\cal Q}_1,{\cal Q}_\infty$.

Notons~:
\begin{eqnarray*}
{\cal V}_0& = & \CC-\RR_{\leq 0},\\
{\cal V}_1&= & \CC-\RR_{\geq 1},\\
{\cal V}_\infty&=&\CC-\RR_{\geq 0}.
\end{eqnarray*}
Une s\'erie localement convergente de ${\cal Q}_0$ d\'etermine de mani\`ere unique 
une fonction analytique dans un voisinage non trivial de $0$ dans ${\cal V}_0$. De m\^{e}me, une s\'erie localement
convergente de 
${\cal Q}_1$ d\'etermine de mani\`ere unique une fonction analytique dans un voisinage non trivial 
de $1$ dans ${\cal V}_1$.
Une s\'erie localement convergente de ${\cal Q}_\infty$ d\'etermine de mani\`ere unique 
une fonction $f$ telle que $f(-1/z)$ soit analytique
dans un voisinage non trivial de $0$ dans ${\cal V}_\infty$. 

\medskip

La fonction $G(z)=z^{\gamma/2}(1-z)^{(\alpha+\beta-\gamma+1)/2}$ est
analytique dans ${\cal V}_0\cap{\cal V}_1$. Dans ce domaine on a~:
\begin{eqnarray*}
G(z) & = & z^{\gamma/2}\sum_{n=0}^\infty\frac{((\alpha+\beta-\gamma+1)/2)_n}{n!}(-1)^nz^n,\quad\mbox{ si }|z|<1\\
& = & z^{\gamma/2}(1+g_1(z)),\\
& = & (1-z)^{(\alpha+\beta-\gamma+1)/2}\sum_{n=0}^\infty\frac{(\gamma/2)_n}{n!}(-1)^n(1-z)^n,\quad\mbox{ si
}|1-z|<1\\
& = & (1-z)^{(\alpha+\beta-\gamma+1)/2}(1+g_2(z)),
\end{eqnarray*}
o\`u $(x)_n=(x-n+1)(x-n+2)\cdots(x-1)x$ est le symbole de Pochammer, et o\`u $g_1\in{\cal M}_0$ et $g_2\in{\cal M}_1$. 
Supposons maintenant que $z\in{\cal V}_0\cap{\cal V}_1$ soit tel que $|z|>1$. On a~:
\begin{eqnarray*}
G(z) & = & e^{\pi{\rm i}\gamma/2}(-z)^{(\alpha+\beta+1)/2}\left(\frac{1}{(-z)}+1\right)^{(\alpha+\beta-\gamma+1)/2}\\
& = & 
\zeta_1(-z)^{(\alpha+\beta+1)/2}\sum_{n=0}^\infty\frac{((\alpha+\beta-\gamma+1)/2)_n}{n!}(-1)^n(-z)^{-n}.
\end{eqnarray*}
Il existe donc $g_3\in{\cal M}_\infty$ tel que 
\[G(z)= \zeta_1(-z)^{(\alpha+\beta+1)/2}(1+g_3(z)).\]
Nous d\'emontrons trois lemmes d\'ecrivant 
les s\'eries de Puiseux de $u_0^2,y_0,y_1,y_2$ au voisinage de $0,1,\infty$.
\subsubsection{S\'eries de Puiseux en $0$.}
\begin{Lemme} Il existe quatre s\'eries localement convergentes $\epsilon_1,\ldots,\epsilon_4\in
{\cal M}_0$ telles que l'on ait,
pour $z\in{\cal V}_0$ et $|z|$ assez petit~:
\begin{eqnarray*}
u_0^2(z) & = & z^{\gamma}(1+\epsilon_1(z)),\\
y_0(z) & = & \frac{\gamma}{2}z^{\gamma-1}(1+\epsilon_2(z)),\\
y_1(z) & = & \frac{\gamma-2}{2}z^{\gamma-1}(1+\epsilon_3(z)),\\
y_2(z) & = & \frac{\gamma}{2}z^{\gamma-1}(1+\epsilon_4(z)).\\
\end{eqnarray*}
\label{proposition:Puiseux_0}\end{Lemme}
\noindent {\bf D\'emonstration}.
Pour $|z|<1$ on a ${}_2F_1(\alpha,\beta;\gamma;z)=1+\kappa_1(z)$ et $\Xi(z)=\gamma^2+\kappa_2(z)$,
avec $\kappa_1,\kappa_2\in{\cal M}_0$. Donc~:
\begin{eqnarray*}
u_0(z) & = & G(z)(1+\kappa_1(z))\\
& = & z^{\gamma/2}(1+g_1(z))(1+\kappa_1(z))\\
& = & z^{\gamma/2}(1+\kappa_3(z)),\\
u_0'(z) & = & \frac{\gamma}{2}z^{(\gamma-2)/2}(1+\kappa_4(z)),
\end{eqnarray*}
pour deux s\'eries localement convergentes $\kappa_3,\kappa_4\in{\cal M}_0$.
On en d\'eduit l'existence de s\'eries localement convergentes $\kappa_5,\ldots,\kappa_{10}\in{\cal M}_0$
telles que~:
\begin{eqnarray*}
u_0(z)^2 & = & z^{\gamma}(1+\kappa_5(z)),\\
y_0(z) & = & u_0(z)u_0'(z)= \frac{\gamma}{2}z^{\gamma-1}(1+\kappa_6(z)),\\
y_1(z) & = & u_0(z)u_0'(z)-u_0(z)^2/z\\
&= & \frac{\gamma}{2}z^{\gamma-1}(1+\kappa_6(z))-z^{\gamma-1}(1+\kappa_7(z))\\
& = & \frac{\gamma-2}{2}z^{\gamma-1}(1+\kappa_8(z)),\\
y_2(z) & = & u_0(z)u_0'(z)-u_0(z)^2/(z-1)\\
& = & \frac{\gamma}{2}z^{\gamma-1}(1+\kappa_6(z))-z^{\gamma}(1+\kappa_9(z))\\
& = & \frac{\gamma}{2}z^{\gamma-1}(1+\kappa_{10}(z)),
\end{eqnarray*}
et nous pouvons poser $\epsilon_1=\kappa_5,\epsilon_2=\kappa_6,\epsilon_3=\kappa_8,\epsilon_4=\kappa_{10}$.

\subsubsection{S\'eries de Puiseux en $1$.}

\begin{Lemme}
Posons
$\theta=\displaystyle{\frac{\Gamma(\gamma)\Gamma(\gamma-\alpha-\beta)}
{\Gamma(\gamma-\alpha)\Gamma(\gamma-\beta)}}$. Il existe quatre s\'eries localement 
convergentes 
$\eta_1,\ldots,\eta_4\in{\cal M}_1$ telles que l'on ait, pour $z\in{\cal V}_1$ et $|1-z|$ assez petit~:
\begin{eqnarray*}
u_0^2(z) & = & \theta^2(1-z)^{1+\alpha+\beta-\gamma}(1+\eta_1(z)),\\
y_0(z) & = & -\frac{\theta^2}{2}(1+\alpha+\beta-\gamma)(1-z)^{\alpha+\beta-\gamma}(1+\eta_2(z)),\\
y_1(z) & = & -\frac{\theta^2}{2}(1+\alpha+\beta-\gamma)(1-z)^{\alpha+\beta-\gamma}(1+\eta_3(z)),\\
y_2(z) & = & -\frac{\theta^2}{2}(-1+\alpha+\beta-\gamma)(1-z)^{\alpha+\beta-\gamma}(1+\eta_4(z)).
\end{eqnarray*}
\label{proposition:Puiseux_1}\end{Lemme}
\noindent {\bf D\'emonstration.}
D'apr\`es \cite{Wittaker:Complex}
14$\cdot$53 p. 290, si $z\in{\cal V}_0\cap{\cal V}_1$, alors on a la relation~:
\begin{eqnarray*}
\lefteqn{\Gamma(\gamma-\alpha)\Gamma(\gamma-\beta)\Gamma(\alpha)\Gamma(\beta){}_2F_1(\alpha,\beta;\gamma;z)=}\\
& = & \Gamma(\gamma)\Gamma(\alpha)\Gamma(\beta)\gamma(\gamma-\alpha-\beta)
{}_2F_1(\alpha,\beta;\alpha+\beta-\gamma+1;1-z)+\\ &  &
\Gamma(\gamma)\Gamma(\gamma-\alpha)
\Gamma(\gamma-\beta)\gamma(\alpha+\beta-\gamma)(1-z)^{\gamma-\alpha-\beta}\times\\
& & {}_2F_1(\gamma-\alpha,\gamma-\beta;\gamma-\alpha-\beta+1;1-z).
\end{eqnarray*}
On a donc~:
\begin{eqnarray}
\lefteqn{{}_2F_1(\alpha,\beta;\gamma;1-z)=}\label{eq:utilotta}\\
& & \frac{\Gamma(\gamma)\Gamma(\gamma-\alpha-\beta)}{\Gamma(\gamma-\alpha)\Gamma(\gamma-\beta)}
{}_2F_1(\alpha,\beta;\alpha+\beta-\gamma+1;z)+\nonumber\\
& & \frac{\Gamma(\gamma)\Gamma(\alpha+\beta-\gamma)}{\Gamma(\alpha)\Gamma(\beta)}z^{\gamma-\alpha-\beta}
{}_2F_1(\gamma-\alpha,\gamma-\beta;\gamma-\alpha-\beta+1;z).\nonumber
\end{eqnarray}
L'hypoth\`ese $\gamma>\alpha+\beta$ implique que le deuxi\`eme terme de (\ref{eq:utilotta}) s'identifie, au voisinage de
$z=0$, \`a une s\'erie localement convergente de ${\cal M}_0$.  Ainsi on obtient, pour $|z|$ assez petit~:
\[{}_2F_1(\alpha,\beta;\gamma;1-z) = \theta(1+\kappa_{11}(z)),\] avec $\kappa_{11}\in{\cal M}_0$.

Pour $|1-z|$ assez petit et $z\in{\cal V}_1$ on a~:
\begin{eqnarray*}
u_0(z) & = & \theta(1-z)^{(\alpha+\beta-\gamma+1)/2}(1+\lambda_{2}(z))\\
u_0'(z) & = & -\frac{\theta}{2}(1-z)^{(\alpha+\beta-\gamma-1)/2}(\alpha+\beta-\gamma+1)(1+\lambda_3(z)),\\
\end{eqnarray*} o\`u $\lambda_2,\lambda_3\in{\cal M}_1$.
Ainsi on obtient les expressions (avec des s\'eries localement convergentes $\lambda_4,\ldots,\lambda_9\in{\cal M}_1$)~:
\begin{eqnarray*}
u_0^2(z) & = & \theta^2(1-z)^{\alpha+\beta-\gamma+1}(1+\lambda_4(z)),\\
y_0(z) & = & -\theta^2\frac{(\alpha+\beta-\gamma+1)}{2}(1-z)^{\alpha+\beta-\gamma}(1+\lambda_5(z))\\
y_1(z) & = &
-\theta^2\frac{(\alpha+\beta-\gamma+1)}{2}(1-z)^{\alpha+\beta-\gamma}(1+\lambda_5(z))-
\\ & & \theta^2(1-z)^{\alpha+\beta-\gamma+1}(1+\lambda_6(z))\\
& = & -\theta^2\frac{(\alpha+\beta-\gamma+1)}{2}(1-z)^{\alpha+\beta-\gamma+1}(1+\lambda_7(z)),\\
y_2(z) & = &
-\theta^2\frac{(\alpha+\beta-\gamma+1)}{2}(1-z)^{\alpha+\beta-\gamma}(1+\lambda_5(z))-\\
& & \theta^2(1-z)^{\alpha+\beta-\gamma-1}(1+\lambda_8(z))\\
& = & -\theta^2\frac{(\alpha+\beta-\gamma-1)}{2}(1-z)^{\alpha+\beta-\gamma+1}(1+\lambda_9(z)),
\end{eqnarray*}
et nous pouvons poser $\eta_1=\lambda_4,\eta_2=\lambda_5,\eta_3=\lambda_7,\eta_4=\lambda_9$.

\subsubsection{S\'eries de Puiseux en $\infty$.}

\begin{Lemme}
Posons
$\omega=\displaystyle{\frac{\Gamma(\gamma)\Gamma(\beta-\alpha)}
{\Gamma(\gamma-\alpha)\Gamma(\beta)}}$ et $\zeta_1=e^{\pi{\rm i}\gamma/2}$. 
Il existe quatre s\'eries localement convergentes
$\mu_1,\ldots,\mu_4\in{\cal M}_{\infty}$ telles
que l'on ait, pour $z\in{\cal V}_\infty$ et $|z|$ assez grand~:
\begin{eqnarray*}
u_0^2(z) & = & (\zeta_1\omega)^2 (-z)^{(1-\alpha+\beta)}(1+\mu_1(z)),\\
y_0(z) & = & \frac{(\zeta_1\omega)^2(\alpha-\beta-1)}{2}(-z)^{\beta-\alpha}(1+\mu_2(z))\\
y_1(z) & = & \frac{(\zeta_1\omega)^2(\alpha-\beta+1)}{2}(-z)^{\beta-\alpha}(1+\mu_3(z))\\
y_2(z) & = & \frac{(\zeta_1\omega)^2(\alpha-\beta+1)}{2}(-z)^{\beta-\alpha}(1+\mu_4(z)).
\end{eqnarray*}
\label{proposition:Puiseux_infty}\end{Lemme}

\noindent {\bf D\'emonstration.} On suppose que $z\in{\cal V}_\infty$. 
D'apr\`es \cite{Wittaker:Complex} 14$\cdot$51 p. 289, On a la formule suivante~:
\begin{eqnarray*}
\lefteqn{\frac{\Gamma(\alpha)\Gamma(\beta)}{\Gamma(\gamma)}{}_2F_1(\alpha,\beta;\gamma;z)=}\\
& = & \frac{\Gamma(\alpha)\Gamma(\beta-\alpha)}{\Gamma(\gamma-\alpha)}(-z)^{-\alpha}
{}_2F_1(\alpha,1-\gamma+\alpha;1-\beta+\alpha;z^{-1})+\\
& & \frac{\Gamma(\beta)\Gamma(\alpha-\beta)}{\Gamma(\gamma-\beta)}(-z)^{-\beta}
{}_2F_1(\beta,1-\gamma+\beta;1-\alpha+\beta;z^{-1}).
\end{eqnarray*}
Donc~:
\begin{eqnarray}
\lefteqn{{}_2F_1(\alpha,\beta;\gamma;z)=}\label{eq:utilotta2}\\
& = & \frac{\Gamma(\gamma)\Gamma(\beta-\alpha)}{\Gamma(\gamma-\alpha)\Gamma(\alpha)}
(-z)^{-\alpha}{}_2F_1(\alpha,1-\gamma+\alpha;1-\beta+\alpha;z^{-1})+\nonumber\\
& & \frac{\Gamma(\gamma)\Gamma(\alpha-\beta)}{\Gamma(\gamma-\beta)\Gamma(\alpha)}(-z)^{-\beta}
{}_2F_1(\beta,1-\gamma+\beta;1-\alpha+\beta;z^{-1}),\nonumber
\end{eqnarray}
Pour $|z|$ assez grand et $z\in{\cal V}_\infty$~:
\begin{eqnarray*}
{}_2F_1(\alpha,\beta;\gamma;z) & = &
\omega(-z)^{-\alpha}( 
{}_2F_1(\alpha,1-\gamma+\alpha;1-\beta+\alpha;z^{-1})+\\
& & (-z)^{\alpha-\beta}\omega_1\omega^{-1}{}_2F_1(\beta,1-\gamma+\beta;1-\alpha+\beta;z^{-1}))\\
& = & \omega(-z)^{-\alpha}(1+\nu_1(z)),
\end{eqnarray*}
avec $\omega=\displaystyle{\frac{\Gamma(\gamma)\Gamma(\beta-\alpha)}{\Gamma(\gamma-\alpha)\Gamma(\alpha)}}$,
$\omega_1=\displaystyle{\frac{\Gamma(\gamma)\Gamma(\alpha-\beta)}{\Gamma(\gamma-\beta)\Gamma(\alpha)}}$, et
$\nu_1\in{\cal M}_\infty$ localement convergente. 

On a~:
\begin{eqnarray*}
u_0(z) & = & \zeta_1(-z)^{(1+\alpha+\beta)/2}(1+g_3(z))\omega (-z)^{-\alpha}(1+\nu_{1}(z))\\
 & = & \zeta_1\omega (-z)^{(1-\alpha+\beta)/2}(1+\nu_{8}(z)),\\
u_0'(z) & = & \frac{\zeta_1\omega}{2}(-z)^{(\beta-\alpha-1)/2}(\alpha-\beta-1)(1+\nu_{9}(z)),\\
\end{eqnarray*} o\`u les $\nu_i$ sont des s\'eries localement convergentes de ${\cal M}_\infty$.
D'ici on obtient~:
\begin{eqnarray*}
u_0^2(z) & = & (\zeta_1\omega)^2 (-z)^{(1-\alpha+\beta)}(1+\nu_{10}(z))\\
y_0(z) & = & u_0(z)u_0'(z)\\
& = &\frac{(\zeta_1\omega)^2(\alpha-\beta-1)}{2}(-z)^{\beta-\alpha}(1+\nu_{11}(z))\\
y_1(z) & = & u_0(z)u_0'(z)+\frac{u_0(z)^2}{(-z)}\\
& = & y_0(z)+(\zeta_1\omega)^2(-z)^{\beta-\alpha}(1+\nu_{12}(z))\\
& = & \frac{(\zeta_1\omega)^2(\alpha-\beta+1)}{2}(-z)^{\beta-\alpha}(1+\nu_{13}(z))\\
y_2(z) & = & u_0(z)u_0'(z)+\frac{u_0(z)^2}{(-z+1)}\\
& = & y_0(z)+(\zeta_1\omega)^2(-z)^{\beta-\alpha}(1+\nu_{14}(z))\\
& = & \frac{(\zeta_1\omega)^2(\alpha-\beta+1)}{2}(-z)^{\beta-\alpha}(1+\nu_{15}(z)).
\end{eqnarray*}
Nous pouvons poser $\mu_1=\nu_{10},\mu_2=\nu_{11},\mu_3=\nu_{13}$ et $\mu_4=\nu_{15}$. 

\medskip

Noter que les nombres 
$\theta$ et $\omega$ sont des coefficients (non nuls) de {\em matrices de connexion} (cf. \cite{Iwasaki:Gauss}).

\section{Ideaux stables.\label{section:ideaux_stables}}
Nous d\'emontrons la proposition suivante.
\begin{Proposition}
Soit ${\cal I}$ un id\'eal premier de ${\cal N}$ non nul et
$D$-stable. Alors ${\cal I}$ contient au moins un des \'el\'ements 
$q,y_0-y_1,y_0-y_2,y_1-y_2$.
\label{proposition:stables}\end{Proposition}

Ainsi, le polyn\^{o}me $\kappa$ de la proposition
\ref{proposition:prop_ramanujan1} est \'egal \`a \[q(y_0-y_1)(y_0-y_2)(y_1-y_2).\]
Pour d\'emontrer la proposition \ref{proposition:stables}, nous devons proc\'eder en plusieurs \'etapes. 
Commen\c{c}ons avec les id\'eaux principaux $D$-stables.

\subsection{Id\'eaux principaux stables.}
On associe \`a $y_i$ le poids $p(y_i)=1$, $i=0,1,2$
et \`a $\tau,q$ le poids $0$. 
Un \'el\'ement de ${\cal N}$ est dit {\em isobare de poids $s$} s'il est somme de mon\^{o}mes
\[c\tau^aq^by_0^{t_0}y_1^{t_1}y_2^{t_2},\] avec $t_0+t_1+t_2=s$ et $c\in\CC^\times$. 
Si $X\in{\cal N}$, nous notons $p(X)$
le plus grand poids d'un mon\^{o}me de $X$. Dans la suite, nous posons aussi 
\[{\cal R}=\CC[y_0,y_1,y_2],\quad{\cal L}=\CC[\tau,q].\] 

On remarque que si $X\in{\cal R}$ 
est non nul et isobare de poids $s$, alors $DX$ est isobare, et
$p(DX)=p(X)+1$. Si $X\not\in{\cal R}$ ou si $X$ n'est pas isobare, on a seulement $p(DX)\leq p(X)+1$
(car les \'el\'ements de poids $0$ de ${\cal N}$ ne sont pas for\c{c}ement tous annul\'es par $D$). D'autre part, pour
$X\in{\cal N}$ non nul, on a clairement
$\deg_q(DX)\leq\deg_q(X)$ et $\deg_\tau(DX)\leq \deg_\tau(X)$.

\begin{Proposition}
Soit ${\cal I}$ un id\'eal premier principal non nul, $D$-stable de 
${\cal N}$. Alors ${\cal I}$ est \'egal \`a l'un des id\'eaux
principaux $(q),(y_0-y_1),(y_0-y_2),(y_1-y_2)$.
\label{lemme:ideaux_principaux_hypergeo}
\end{Proposition}
Cette proposition g\'en\'eralise le lemme 5.2 p. 161 de \cite{Nesterenko:Introduction10}.

\medskip

\noindent {\bf D\'emonstration}. 
En utilisant les formules explicites du lemme \ref{lemme:relations_explicites}, on remarque pour commencer que 
\begin{eqnarray*}
Dq & = & wq,\\
D(y_0-y_1)& = & (y_0-y_1)(y_0+y_1),\\
D(y_0-y_2)& = & (y_0-y_2)(y_0+y_2),\\
D(y_1-y_2)& = &(y_1-y_2)(y_1+y_2).\end{eqnarray*}
Donc les id\'eaux principaux $(q),(y_0-y_1),(y_0-y_2),(y_1-y_2)$
de ${\cal N}$, qui sont premiers, sont aussi $D$-stables.

Soit $P\in{\cal N}$ un polyn\^{o}me non nul et irr\'eductible tel que 
\begin{equation}
DP=FP\label{eq:principalite}
\end{equation} avec $F\in{\cal N}$. Nous devons d\'emontrer que $P$ est proportionnel \`a l'un des
polyn\^{o}mes 
$q,y_0-y_1,y_0-y_2,y_1-y_2$, avec une constante de proportionnalit\'e dans $\CC^\times$.

Nous d\'eterminons une expression plus pr\'ecise de $F$.
Comme $p(DP)=p(F)+p(P)$, on a que $p(F)\leq 1$. Comme
de plus $\deg_q(DP)=\deg_q(F)+\deg_q(P)$ et 
$\deg_\tau(DP)=\deg_\tau(F)+\deg_\tau(P)$, on a que $\deg_q(F)=\deg_\tau(F)=0$ et donc
\begin{equation}F=\lambda_0y_0+\lambda_1y_1+\lambda_2y_2+\lambda_3,\label{eq:forme_f}\end{equation} avec
$\lambda_0,\lambda_1,\lambda_2,\lambda_3\in\CC$. On a aussi $F\not=0$ car 
$P$ \'etant irr\'eductible, il est non constant.

Notre premier objectif est de montrer que dans (\ref{eq:forme_f}) on a 
$\lambda_i\in\QQ$, $i=0,\ldots,3$. 

\begin{Lemme} Sous l'hypoth\`ese (\ref{eq:principalite}), on a 
$\lambda_0,\lambda_1,\lambda_2\in\QQ$ dans (\ref{eq:forme_f}).\label{lemme:123}\end{Lemme}

\noindent {\bf D\'emonstration.}
Ecrivons~:
\begin{equation}
P=P_h+P_{h+1}+\cdots+P_k,\label{eq:forme}\end{equation} avec $P_i$ isobare de poids $i$, $P_h$ non nul isobare de poids minimal
$h\leq k$, et $P_k$ non nul isobare de poids maximal
$k$. Ecrivons aussi~:
\[P_k=\sum_{j=1}^tf_jV_j,\]
avec $V_1,\ldots,V_t$ isobares dans ${\cal R}$ de poids $k$, et
$f_j\in{\cal L}$. 
On a que $DP_k$ est somme d'un \'el\'ement isobare de ${\cal N}$ de poids $k+1$ et d'un \'el\'ement isobare de ${\cal N}$
de poids $k$~; de plus $P$ satisfait (\ref{eq:principalite}). Nous obtenons, en tenant compte des termes de poids $k$ et
$k+1$ de
$FP$~:
\begin{eqnarray*}
DP_k & = & \sum_{j=1}^tf_jDV_j+\sum_{j=1}^t(Df_j)V_j\\
& = & \lambda_3P_k+P_{k-1}\sum_{i=0}^2\lambda_iy_i+P_k\sum_{i=0}^2\lambda_i y_i.
\end{eqnarray*}
En comparant les poids des termes
isobares de poids $k+1$, on voit que~:
\begin{eqnarray*}
\sum_{j=1}^tf_jDV_j & = & \frac{\partial P_k}{\partial y_0}Dy_0+\frac{\partial P_k}{\partial y_1}Dy_1+
\frac{\partial P_k}{\partial y_2}Dy_2\\
& = & P_k\sum_{i=0}^2\lambda_i y_i.
\end{eqnarray*}
Soit $D'$ la d\'erivation de ${\cal N}$ d\'efinie par $D'=D$ sur ${\cal R}$ et $D'$ nulle sur ${\cal L}$.
On a donc d\'emontr\'e que 
\begin{equation}
D'P_k=P_k\sum_{i=0}^2\lambda_i y_i.\label{eq:derivation_prime}\end{equation}
Nous \'etudions l'influence de la r\'elation (\ref{eq:derivation_prime}) sur des s\'eries de Puiseux formel\-les.
Nous allons d\'eduire de (\ref{eq:derivation_prime}) et des lemmes
\ref{proposition:Puiseux_0}, \ref{proposition:Puiseux_1} et \ref{proposition:Puiseux_infty}
les trois relations suivantes~:
\begin{small}\begin{eqnarray}
\lambda_0\gamma+\lambda_1(\gamma-2)+\lambda_2\gamma & \in & \QQ,\label{eq:0}\\
\lambda_0(\alpha+\beta-\gamma+1)+\lambda_1(\alpha+\beta-\gamma+1)+\lambda_2(\alpha+\beta-\gamma-1) & \in & \QQ,\label{eq:1}\\
\lambda_0(\alpha-\beta-1)+\lambda_1(\alpha-\beta+1)+\lambda_2(\alpha-\beta+1) & \in & \QQ.\label{eq:infty}
\end{eqnarray}\end{small}
D\'emontrons l'\'egalit\'e (\ref{eq:0}). 
D'apr\`es le lemme \ref{proposition:Puiseux_0} et l'ind\'ependance alg\'ebri\-que sur $\CC$ de $y_0,y_1,y_2$,
il existe un homomorphisme injectif naturel~:
\[\Phi_0:{\cal R}\rightarrow{\cal Q}_0,\] qui est d\'efini en associant \`a tout \'el\'ement 
de ${\cal R}$ la s\'erie de Puiseux donn\'ee \`a partir du lemme \ref{proposition:Puiseux_0}.

Soit ${\cal R}_0$ le sous-anneau de ${\cal Q}_0$ image de $\Phi_0$~; on a un isomorphisme d'anneaux $\Phi_0:{\cal
R}\cong{\cal R}_0$. 
Soient $U,V$ deux ind\'etermin\'ees et posons~:
\[{\cal N}_0=\CC[U,V,\Phi_0(y_0),\Phi_0(y_1),\Phi_0(y_2)].\]
Autrement dit,
${\cal N}_0$ est le sous-anneau des s\'eries de Puiseux formelles en
les puissances enti\`eres de $z^{1/p}$ \`a coefficients des polyn\^{o}mes en
$U,V$, qui est engendr\'e par les s\'eries de Puiseux 
de $y_0,y_1,y_2$ donn\'ees par le lemme \ref{proposition:Puiseux_0}.

L'isomorphisme $\Phi_0:{\cal
R}\cong{\cal R}_0$ se prolonge en un isomorphisme $\Phi_0:{\cal N}\rightarrow{\cal N}_0$
en posant $\Phi_0(\tau)=U$ et $\Phi_0(q)=V$, car $\tau,q,y_0,y_1,y_2$ sont alg\'ebriquement ind\'ependantes sur $\CC$,
comme nous l'avons rappel\'e dans l'introduction.

La d\'erivation $D:{\cal R}\rightarrow{\cal R}$ d\'ecrite dans le lemme \ref{lemme:relations_explicites}
d\'etermine une d\'erivation $D_0:{\cal R}_0\rightarrow{\cal R}_0$ d\'efinie sur les s\'eries de Puiseux formelles
de la mani\`ere suivante. Notons $\sigma_0$ la s\'erie de Puiseux formelle de $u_0^2$ donn\'ee par le lemme 
\ref{proposition:Puiseux_0}, soit
$\varphi=\sum_mf_mz^{m/p}$ une s\'erie de Puiseux de
${\cal R}_0$ (ou plus g\'en\'eralement de ${\cal Q}_0$). On pose~:
\begin{equation}
D_0(\varphi) = \sigma_0\sum_m(m/p)f_mz^{(m/p)-1}.\label{eq:def_deri_0}\end{equation}
La d\'erivation $D_0:{\cal R}_0\rightarrow{\cal R}_0$ induit une d\'erivation 
$D_0:{\cal N}_0\rightarrow{\cal N}_0$ en posant $D_0(U)=D_0(V)=0$~:
elle se calcule sur les s\'eries formelles en utilisant (\ref{eq:def_deri_0}), et en permettant que 
$f_m$ appartienne aussi \`a
$\CC[U,V]$. On v\'erifie facilement que si $s\in{\cal N}$, alors $\Phi_0(D'(s))=D_0(\Phi_0(s))$~:
on a ainsi construit un isomorphisme d'anneaux diff\'erentiels~:
\begin{equation}\Phi_0:({\cal N},D')\rightarrow({\cal N}_0,D_0).\label{eq:iso_diff_0}\end{equation}

Posons~: \[\Phi_0(P_k) =\sum_{m\geq m_0}a_mz^{m/p}\] avec $a_m\in\CC[U,V]$ et $a_{m_0}$ non nul. On a,
d'apr\`es le lemme \ref{proposition:Puiseux_0}~:
\begin{eqnarray*}
D_0\left(\sum_{m\geq m_0}a_mz^{m/p}\right) & = & \sigma_0\sum_{m\geq m_0}(m/p)a_mz^{(m/p)-1}\\
& = & z^\gamma(1+\Phi_0(\epsilon_1(z)))\sum_{m\geq m_0}(m/p)a_mz^{(m/p)-1}\\
& = & (1+\epsilon_5(z))(m_0/p)a_{m_0}z^{\gamma+(m_0/p)-1},
\end{eqnarray*} avec $\epsilon_5\in{\cal M}_0{\cal N}_0\subset{\cal N}_0$.
L'\'egalit\'e
(\ref{eq:derivation_prime}), l'isomorphisme (\ref{eq:iso_diff_0}) et 
le lemme \ref{proposition:Puiseux_0} impliquent~:
\begin{eqnarray*}
D_0\left(\sum_{m\geq m_0}a_mz^{m/p}\right) & = & \sum_i\lambda_i\Phi_0(y_i)\sum_{m\geq m_0}a_mz^{m/p}\\ 
& = & \frac{1}{2}(\lambda_0\gamma+\lambda_1(\gamma-2)+\lambda_2\gamma+\epsilon_6(z))\times\\
& & z^{\gamma-1}\sum_{m\geq
m_0}a_mz^{m/p},
\end{eqnarray*}
o\`u $\epsilon_6(z)$ est un autre \'el\'ement de ${\cal M}_0{\cal N}_0$.
On en d\'eduit l'identit\'e formelle~:
\begin{eqnarray*}
z^\gamma (m_0/p)a_{m_0}z^{(m_0/p)-1} & = & \frac{1}{2}(\lambda_0\gamma+\lambda_1(\gamma-2)+\lambda_2\gamma)z^{\gamma-1}
a_{m_0}z^{m_0/p},
\end{eqnarray*}
d'o\`u la relation (\ref{eq:0}), car $a_{m_0}\not=0$ et $m_0\in\ZZ$.

Les r\'elations (\ref{eq:1}) et (\ref{eq:infty}) se d\'emontrent de mani\`ere similaire.
Soit ${\cal R}_1$ le sous-anneau de ${\cal Q}_1$ engendr\'e par les s\'eries de Puiseux formelles de $y_0,y_1,y_2$
du lemme \ref{proposition:Puiseux_1}, notons \[{\cal N}_1=\CC[U,V,\Phi_1(y_0),\Phi_1(y_1),\Phi_1(y_2)],\]
soit $\sigma_1\in{\cal Q}_1$ la s\'erie de Puiseux de $u_0^2$ au voisinage de $1$ donn\'ee par le lemme
\ref{proposition:Puiseux_1}, soit $D_1$ la d\'erivation de ${\cal N}_1$ d\'efinie de la mani\`ere suivante. Pour
$\varphi=\sum_{m}f_m(1-z)^{m/p}$ un \'el\'ement de ${\cal N}_1$, on pose~:
\[D_1(\varphi)=-\sigma_1\sum_{m}(m/p)f_m(1-z)^{(m/p)-1}.\] On a un isomorphisme d'anneaux diff\'erentiels
$\Phi_1:({\cal N},D')\cong({\cal N}_1,D_1)$. 

Posons~:
\[\Phi_1(P_k)=\sum_{m\geq m_1}b_m(1-z)^{m/p},\] avec $b_m\in\CC[U,V]$ et $b_{m_1}\not=0$. L'\'egalit\'e
(\ref{eq:derivation_prime}) implique~:\[
D_1\left(\sum_{m\geq m_1}b_m(1-z)^{m/p}\right)=\sum_i\lambda_i\Phi_1(y_i)\sum_{m\geq
m_1}b_m(1-z)^{m/p}\] d'o\`u l'on d\'eduit, en suivant les m\^{e}mes techniques que nous avons employ\'e pour obtenir la relation
(\ref{eq:0}), et en appliquant les formules du lemme \ref{proposition:Puiseux_1}~:
\begin{eqnarray*}
\lefteqn{-\theta^2(1-z)^{1+\alpha+\beta-\gamma} (m_1/p)b_{m_1}(1-z)^{(m_1/p)-1}=
-\frac{\theta^2}{2}(1-z)^{\alpha+\beta-\gamma}\times}\\
& &(\lambda_0(1+\alpha+\beta-\gamma)+\lambda_1(1+\alpha+\beta-\gamma)+\lambda_2(-1+\alpha+\beta-\gamma))\times\\
& & b_{m_1}(1-z)^{m_1/p},
\end{eqnarray*}
d'o\`u\begin{small}
\[\frac{m_1}{p}=(1/2)(\lambda_0(1+\alpha+\beta-\gamma)+\lambda_1(1+\alpha+\beta-\gamma)+\lambda_2(-1+\alpha+\beta-\gamma)),\]
\end{small}
car $\theta,b_{m_1}\not=0$, ce qui implique la r\'elation (\ref{eq:1}).

Pour d\'emontrer la r\'elation (\ref{eq:infty}) on introduit le sous-anneau ${\cal R}_\infty$ de ${\cal Q}_\infty$
engendr\'e par les s\'eries de Puiseux formelles de $y_0,y_1,y_2$
du lemme \ref{proposition:Puiseux_infty}, et on note \[{\cal
N}_\infty=\CC[U,V,\Phi_\infty(y_0),\Phi_\infty(y_1),\Phi_\infty(y_2)],\] 
puis on d\'esigne par $\sigma_\infty\in{\cal Q}_\infty$
la s\'erie de Puiseux de $u_0^2$ au voisinage de $\infty$ donn\'ee par  le lemme
\ref{proposition:Puiseux_infty}. Nous construisons, comme nous l'avons fait
pour $D_0$ et $D_1$, une d\'erivation $D_\infty:{\cal N}_\infty\rightarrow{\cal N}_\infty$ telle qu'on ait 
un isomorphisme $\Phi_\infty:({\cal N},D')\cong({\cal N}_\infty,D_\infty)$.
Explicitement, soit $\varphi=\sum_{m}f_m(-z)^{-m/p}$ un \'el\'ement de ${\cal N}_\infty$~; on pose~:
\[D_\infty(\varphi)=-\sigma_\infty\sum_{m}(m/p)f_m(-z)^{-(m/p)-1}.\] Ecrivons~: \[\Phi_\infty(P_k)=\sum_{m\geq
m_\infty}f_m(-z)^{-m/p}.\] L'\'egalit\'e
(\ref{eq:derivation_prime}) \'equivaut \`a~:\[
D_\infty\left(\sum_{m\geq m_\infty}f_m(-z)^{-m/p}\right)=\sum_i\lambda_i\Phi_\infty(y_i)\sum_{m\geq m_\infty}f_m(-z)^{-m/p}\]
d'o\`u l'on d\'eduit, en appliquant  les formules du lemme \ref{proposition:Puiseux_infty}~:
\begin{eqnarray*}
\lefteqn{-(\zeta_1\omega)^2(-z)^{1-\alpha+\beta}(m_\infty/p) f_{m_\infty}(-z)^{(m_\infty/p)-1} =}\\
& = &
\frac{(\zeta_1\omega)^2}{2}(\lambda_0(\alpha-\beta-1)+\lambda_1(\alpha-\beta+1)+\lambda_2(\alpha-\beta+1))
(-z)^{\beta-\alpha}\times\\
& & f_{m_\infty}(-z)^{m_\infty/p},
\end{eqnarray*}
ce qui implique 
\[\frac{m_\infty}{p}=-(1/2)(\lambda_0(\alpha-\beta-1)+\lambda_1(\alpha-\beta+1)+\lambda_2(\alpha-\beta+1)),\]
car $f_{m_\infty}\not=0$ et $\zeta_1\omega\not=0$.
Nous avons obtenu aussi la relation (\ref{eq:infty}).

Posons~:
\[M:=\left(\begin{array}{ccc}
\gamma & \gamma-2 & \gamma \\ \alpha+\beta-\gamma+1 & \alpha+\beta-\gamma+1 & \alpha+\beta-\gamma-1 \\
\alpha-\beta-1 & \alpha-\beta+1 & \alpha-\beta+1
\end{array}\right).\]
Les r\'elations (\ref{eq:0}), (\ref{eq:1}), (\ref{eq:infty}) que nous avons
d\'emontr\'e, impliquent qu'il existe des nombres rationnels $\mu_0,\mu_1,\mu_2$ tels que (le symbole ${}^t$ d\'esigne la
transpos\'ee de matrice)~:
\[{}^t(\mu_0,\mu_1,\mu_2)=M\cdot{}^t(\lambda_0,\lambda_1,\lambda_2).\] La matrice $M$ a d\'eterminant $8\alpha\not=0$
et a tous ses coefficients rationnels.
Soit $N$ son inverse. Alors~:
\[{}^t(\lambda_0,\lambda_1,\lambda_2)=N\cdot{}^t(\mu_0,\mu_1,\mu_2)\in\QQ^3,\] nous avons donc d\'emontr\'e que 
$\lambda_0,\lambda_1,\lambda_2\in\QQ$.

\begin{Lemme}
sous l'hypoth\`ese (\ref{eq:principalite}), on a $\lambda_3\in\QQ$
dans (\ref{eq:forme_f}).
\label{lemme:1}\end{Lemme}
\noindent {\bf D\'emonstration.}
Consid\'erons \`a nouveau l'expression (\ref{eq:forme}), o\`u nous avions suppos\'e
$P_h\not=0$ et \'ecrivons~:
\[P_h=\sum_{j=1}^td_jU_j,\] avec $U_1,\ldots,U_t$ \'el\'ements isobares de 
${\cal R}$ de poids
$h$, et $d_j\in{\cal L}$. Si $h=0$, alors $P_0,DP_0\in{\cal L}$. Si $h\not=0$, alors on a
que $DP_h$ est somme d'un \'el\'ement de poids $h$ et un \'el\'ement de poids $h+1$. En regardant les 
termes de poids $h$ et $h+1$ dans l'expression de $FP$ comme somme d'\'el\'ements isobares, nous trouvons
que~:
\begin{eqnarray*}
DP_h & = & \sum_{j=1}^t(Dd_j)U_j+\sum_{j=1}^td_jDU_j\\
& = & \lambda_3P_h+\lambda_3P_{h+1}+P_h\sum_{i=0}^2\lambda_iy_i.
\end{eqnarray*}
Si $h=0$, alors nous trouvons la relation plus simple $DP_0=\lambda_3P_0$. Dans tous les cas, en comparant 
les termes isobares de poids $h$~:
\begin{eqnarray*}
\sum_{j=1}^t(Dd_j)U_j & = & w\left(\frac{\partial P_h}{\partial\tau}+\frac{\partial P_h}{\partial q}q\right)\\
& = & \lambda_3P_h.
\end{eqnarray*}
Soit $H$ la d\'erivation de ${\cal N}$ d\'efinie par $H({\cal R})=0$ et $H=D$ sur ${\cal L}$. 
Nous avons~:
\[H(P_h)=\lambda_3P_h.\]
Ecrivons $P_h$ explicitement~:
\begin{eqnarray*}
P_h & = & \sum_{a,b\geq 0}c_{a,b}\tau^aq^b\\
& =& \sum_{a,b>0}c_{a,b}\tau^aq^b+\sum_{b>0}d_bq^b+\sum_{a>0}e_a\tau^a+f_0,\quad c_{a,b},d_b,e_a,f_0\in{\cal R}.
\end{eqnarray*}
On a~:
\begin{eqnarray*}
w^{-1}H(P_h) & = & 
\sum_{a,b}c_{a,b}(a\tau^{a-1}q^b+b\tau^aq^b)\\ & = &
\sum_{a,b>0}(bc_{a,b}+(a+1)c_{a+1,b})\tau^aq^b+\sum_{b>0}(bd_b+c_{1,b})q^b+\\
& & \sum_{a>0}(a+1)e_{a+1}\tau^a+e_1.
\end{eqnarray*}
Soit $a_0$ maximal avec $c_{a_0,b_0}\not=0$ pour quelques $b_0\geq 0$.
Si $a_0=0$, alors $P_h\in{\cal R}[q]$. Donc \[P_h=U_0+U_1q+\cdots+U_rq^r\] avec $U_j\in{\cal R}$.
Si $r=0$ alors $P_h\in{\cal R}$ et $\lambda_3=0$. Sinon,
$r>0$, $U_r\not=0$, et
$w^{-1}H(P_h)=k_1q+\cdots+rk_rq^r$, d'o\`u $rk_r=w\lambda_3k_r$ ce qui implique $\lambda_3=rw^{-1}\in\QQ$.

Si $a_0>0$ alors $b_0>0$, car sinon, $\deg_\tau(H(P_h))<\deg_\tau(P_h)$ et nous ne pouvons pas avoir
$H(P_h)=\lambda_3P_h$. Dans ce cas, $c_{a_0+1,b_0}=0$.
Donc~:
\begin{eqnarray*}
w^{-1}b_0c_{a_0,b_0} & = &  
b_0c_{a_0,b_0}+(a_0+1)c_{a_0+1,b_0}\\
& = & \lambda_3c_{a_0,b_0},
\end{eqnarray*}
d'o\`u $\lambda_3=w^{-1}b_0$ et $\lambda_3\in\QQ$.

\subsubsection{Fin de la preuve de la proposition \ref{lemme:ideaux_principaux_hypergeo}.}
Soit $P$ un polyn\^{o}me irr\'edu\-ctible de ${\cal N}$ satisfaisant (\ref{eq:principalite}), avec $F$ comme dans
(\ref{eq:forme_f}), c'est-\`a-dire, tel que
$DP=(\lambda_0y_0+\lambda_1y_1+\lambda_2y_2+\lambda_3)P$. A pr\'esent, nous avons d\'emontr\'e que
$\lambda_i\in\QQ$ pour $i=0,\ldots,3$. Donc
\[\displaystyle{\frac{DP}{P}},\displaystyle{\frac{D(y_0-y_1)}{(y_0-y_1)}},\displaystyle{\frac{D(y_0-y_2)}{(y_0-y_2)}},
\displaystyle{\frac{D(y_1-y_2)}{(y_1-y_2)}},\displaystyle{\frac{Dq}{q}}\] sont $\QQ$-lin\'eairement d\'ependants, et il existe des
entiers rationnels non tous nuls $\alpha_1,\ldots,\alpha_5$ tels que~:
\[\alpha_1\frac{DP}{P}+\alpha_2\frac{D(y_0-y_1)}{(y_0-y_1)}+\alpha_3\frac{D(y_0-y_2)}{(y_0-y_2)}+\alpha_4
\frac{D(y_1-y_2)}{(y_1-y_2)}+\alpha_5\frac{Dq}{q}=0.\]
Ainsi $D(P^{\alpha_1}(y_0-y_1)^{\alpha_2}(y_0-y_2)^{\alpha_3}(y_1-y_2)^{\alpha_4}q^{\alpha_5})=0$, d'o\`u
\[P^{\alpha_1}(y_0-y_1)^{\alpha_2}(y_0-y_2)^{\alpha_3}(y_1-y_2)^{\alpha_4}q^{\alpha_5}\in\CC^\times.\] Comme $P$ est
irr\'eductible,
ceci implique que $P$ est proportionnel \`a l'un des polyn\^{o}mes $y_0-y_1,y_0-y_2,y_1-y_2,q$.
La d\'emonstration du lemme est maintenant termin\'ee.

\subsection{Id\'eaux non principaux stables.}

Nous introduisons ici les crochets de Rankin, qui sont des op\'erateurs diff\'erentiels
d\'efinis uniquement sur des polyn\^{o}mes isobares de ${\cal R}$, mais pas sur des polyn\^{o}mes quelconque de ${\cal N}$~:
il permettent d'\'eviter les arguments du paragraphe 5 pp. 162-165 de \cite{Nesterenko:Introduction10}. 

Le crochet de Rankin $[U,V]$ de deux polyn\^{o}mes isobares $U,V$ de ${\cal R}$ est 
d\'efini par \[[U,V]=\displaystyle{p(U) U D(V)-p(V) V D(U)}.\]
Voici les propri\'et\'es principales du crochet.

\begin{Lemme} Soient $X,Y,P$ des \'el\'ements isobares de ${\cal R}$.
\begin{itemize}
\item Si $X\in{\cal R}$ est de poids $x$ et $Y$ est de poids $y$, alors $[X,Y]=-[Y,X]$ est isobare de poids $x+y+1$.
\item Si $X,Y$ ont m\^{e}me poids, alors l'application $d_P(X):=[X,P]$
satisfait $d_P(X+Y)=d_P(X)+d_P(Y)$. 
\item On a $d_P(XY)=d_P(X)Y+Xd_P(Y)$.
\item Si ${\cal I}$ est un id\'eal de ${\cal R}$ qui est $D$-stable, et $X\in{\cal I}$, alors $d_P(X)\in{\cal I}$.
\end{itemize}\label{lemme:proprietes_crochet}\end{Lemme}
Les premi\`eres deux propri\'et\'es sont triviales. D\'emontrons la troisi\`eme propri\'et\'e. On a~:
\begin{eqnarray*}
d_P(XY)=[XY,P] & = & p(XY)XYDP-p(P)PD(XY)\\
& = & (p(X)+p(Y))XYDP-p(P)P((DX)Y+XDY)\\
& = & (p(X)XDP-p(P)PDX)Y+\\ & & X(p(Y)YDP-p(P)PDY)\\
& = & d_P(X)Y+Xd_P(Y).
\end{eqnarray*}
D\'emontrons la quatri\`eme propri\'et\'e. Si $X\in{\cal I}$ et $P$ est isobare, alors $DX\in{\cal I}$ car ${\cal I}$ est $D$-stable. Donc
$XDP,PDX\in{\cal I}$, et $d_P(X)=p(X)XDP-p(P)PDX\in{\cal I}$.

\subsubsection{Id\'eaux stables isobares.}

Nous d\'emontrons maintenant le lemme suivant.
\begin{Lemme}
Tout id\'eal premier isobare de ${\cal R}$, $D$-stable, contient au moins un des \'el\'ements 
$y_0-y_1,y_0-y_2,y_1-y_2$.\label{lemme:une_variable_hypergeo}\end{Lemme}
{\bf D\'emonstration}. Soit ${\cal I}$ un id\'eal comme dans les hypoth\`eses du lemme. 
Si ${\cal I}$ est principal alors ${\cal I}{\cal N}$ est premier, principal et $D$-stable,
et au moins un des polyn\^{o}mes $y_0-y_1,y_0-y_2,y_1-y_2$
appartient \`a ${\cal I}$, d'apr\`es la proposition \ref{lemme:ideaux_principaux_hypergeo}. Supposons donc que ${\cal I}$
ne soit pas principal: nous avons deux cas.

\medskip

\noindent {\bf (1)}. Supposons que ${\cal I}\cap(\CC[y_0]\cup\CC[y_1]\cup\CC[y_2])\not=\{0\}$.
Supposons par exemple que ${\cal I}\cap\CC[y_0]\not=(0)$. Comme ${\cal I}$ est premier, il existe $v\in\CC$ tel que
$y_0-v\in{\cal I}$~: comme ${\cal I}$ est isobare, on a $v=0$. 

Donc $y_0\in{\cal I}$ et 
$Dy_0\in{\cal I}$. Ainsi le polyn\^{o}me $H$ obtenu en substituant $y_0\mapsto 0$ dans le polyn\^{o}me $Dy_0$
appartient aussi \`a ${\cal I}$. Explicitement~:
\[H=-\frac{1}{4}(ay_1^2+by_2^2+c(y_1-y_2)^2)\in{\cal I}.\] Aussi le polyn\^{o}me $DH$ appartient \`a
${\cal I}$, et le polyn\^{o}me $K$ obtenu en substituant $y_0\mapsto 0$ dans $DH$ appartient \`a ${\cal I}$.
On a~:
\begin{eqnarray*}
K & = & \frac{1}{8}(a^2y_1^3+b(b-4)y_2^3-c(y_1-y_2)^2(4y_1+y_2(4-b))+\\
& & ay_1((c-4)y_1^2+(b-2c)y_1y_2+(b+c)y_2^2))\\
& & \in{\cal I}.
\end{eqnarray*}
Les r\'esultants $R_1=$R\'es${}_{y_1}(H,K)$ et $R_2=$R\'es${}_{y_2}(H,K)$ appartiennent \`a ${\cal I}$.
Explicitement, on trouve $R_1=\displaystyle{-\frac{1}{256}\kappa y_2^6}$ et 
$R_2=\displaystyle{-\frac{1}{256}\kappa y_1^6}$, avec
\[\eta=(a+b)(a+c)(b+c)(ab+bc+ac).\] 
On v\'erifie que $\eta\not=0$ (voir le lemme
\ref{lemme:elementaire}, \`a la fin de ce paragraphe), donc si $y_0\in{\cal I}$, alors aussi
$y_1,y_2\in{\cal I}$. On parvient au m\^{e}me r\'esultat si on suppose que $y_1\in{\cal I}$ ou $y_2\in{\cal I}$. Pour faire les
calculs on peut exploiter les sym\'etries \'evidentes du syst\`eme diff\'erentiel. Ainsi, on a que si ${\cal I}\cap
(\CC[y_0]\cup\CC[y_1]\cup\CC[y_2])\not=\{0\}$, alors $y_0,y_1,y_2\in{\cal I}$, et dans ce cas on obtient la
conclusion du lemme, car on trouve m\^{e}me plus, \`a savoir 
que ${\cal I}\supset(y_0-y_1,y_0-y_2,y_1-y_2)= (y_0-y_1,y_0-y_2)$.

\medskip

\noindent {\bf (2).}
Supposons que ${\cal I}\cap(\CC[y_0]\cup\CC[y_1]\cup\CC[y_2])=\{0\}$.
Alors ${\cal I}$ contient un polyn\^{o}me
isobare $P$ tel que $\displaystyle{\frac{\partial P}{\partial y_1}\not=0}$ et 
$\displaystyle{\frac{\partial P}{\partial y_2}=0}$. En d'autres termes, ${\cal I}\cap\CC[y_0,y_1]\not=(0)$.
Nous pouvons aussi supposer que $P$ ait degr\'e minimal en $y_1$~: donc $\displaystyle{\frac{\partial P}{\partial
y_1}}\not\in{\cal I}$. Nous calculons le crochet de Rankin $[P,y_0]$, qui est 
un polyn\^{o}me non nul appartenant \`a ${\cal I}$ (car $P$ n'est pas une puissance de $y_0$). On obtient~:
\begin{eqnarray*}
[P,y_0] & = & \frac{\partial P}{\partial y_0}[y_0,y_0]+\frac{\partial P}{\partial y_1}[y_1,y_0]+
\frac{\partial P}{\partial y_2}[y_2,y_0]\\
& = & \frac{\partial P}{\partial y_1}[y_1,y_0]\in{\cal I},
\end{eqnarray*}
car $[y_0,y_0]=0$ et $\displaystyle{\frac{\partial P}{\partial y_2}=0}$. Comme $\displaystyle{\frac{\partial P}{\partial
y_1}}\not\in{\cal I}$, on a $[y_1,y_0]\in{\cal I}$.

De m\^{e}me, ${\cal I}\cap\CC[y_0,y_2]\not=(0)$. Soit $Q$ un polyn\^{o}me isobare de ${\cal I}$
tel que $\displaystyle{\frac{\partial Q}{\partial y_2}\not=0}$ et 
$\displaystyle{\frac{\partial Q}{\partial y_1}=0}$. Nous pouvons aussi supposer que $Q$ soit tel
que $\displaystyle{\frac{\partial Q}{\partial y_2}}\not\in{\cal I}$.
En suivant les m\^{e}mes arguments que ci-dessus, on trouve que $[y_2,y_0]\in{\cal I}$. Donc $[y_1,y_0],[y_2,y_0]\in{\cal I}$,
et $[y_1,y_0]-[y_2,y_0]\in{\cal I}$. Mais~:
\begin{eqnarray*}
[y_1,y_0]-[y_2,y_0] & = & [y_1-y_2,y_0]\\
& = & (y_1-y_2)y_0.
\end{eqnarray*}
comme $y_0\not\in{\cal I}$, on obtient $y_1-y_2\in{\cal I}$, et la d\'emonstration du lemme est termin\'ee.

\medskip

\noindent {\bf Remarque.} Les sym\'etries du syst\`eme diff\'erentiel permettent en fait de d\'emontrer que
si ${\cal I}$ est un id\'eal non principal satisfaisant les hypoth\`eses du lemme \ref{lemme:une_variable_hypergeo}, 
alors ${\cal
I}\supset(y_0-y_1,y_0-y_2)$.

\medskip

Afin de d\'emontrer la proposition \ref{proposition:stables}, nous devons 
encore \'etablir un lemme \'el\'ementaire. Soit ${\cal I}$ un id\'eal de ${\cal R}$~; 
on note $\tilde{{\cal I}}$ l'id\'eal engendr\'e par
les poly\-n\^{o}mes isobares de ${\cal I}$.
\begin{Lemme}
Nous avons les propri\'et\'es suivantes.
\begin{enumerate}
\item Si ${\cal I}$ est premier, alors $\tilde{{\cal I}}$ est premier.
\item Si ${\cal I}$ est $D$-stable, alors $\tilde{{\cal I}}$ est $D$-stable.
\item Si ${\cal I}$ n'est pas principal, alors $\tilde{{\cal I}}\not=(0)$.
\end{enumerate}
\label{lemme:nuovo}\end{Lemme}
\noindent {\bf D\'emonstration.} 1. C'est bien connu. 

\noindent 2. Soit $U\in\tilde{{\cal I}}$~; alors $U=U_0+\cdots+U_k$ pour certains
polyn\^{o}mes isobares $U_0,\ldots,U_k\in{\cal I}$. Pour tout $j$ on a $DU_j\in{\cal I}$ car ${\cal I}$ est $D$-stable.
De plus $DU_j$ est isobare, donc $DU=DU_0+\cdots+DU_k\in\tilde{{\cal I}}$.

\noindent 3. Introduisons une inconnue $y_3$ de poids $1$, et notons, pour tout 
polyn\^{o}me non nul $F\in{\cal R}$~:
\[{}^hF=y_3^{p}F\left(\frac{y_0}{y_3},\frac{y_1}{y_3},\frac{y_2}{y_3}\right),\] o\`u $p$ est le plus grand poids d'un mon\^{o}me 
non nul de $P$.

Il faut montrer que ${\cal I}$ contient un polyn\^{o}me isobare non nul. Par hypoth\`ese il existe 
$U,V\in{\cal I}$ premiers entre eux~; donc ${}^hU$ et ${}^hV$ sont premiers entre eux. 

Le r\'esultant
$R=\mbox{R\'es}_{y_3}({}^hU,{}^hV)\in{\cal R}$ est non nul et on v\'erifie qu'il est isobare. De plus $R\in{\cal I}$, car 
il existe $A,B\in{\cal R}[y_3]$ tels que $R=A({}^hU)+B({}^hV)$ d'o\`u, en substituant $y_3=1$, $R=\tilde{A}U+\tilde{B}V$ pour
$\tilde{A},\tilde{B}\in{\cal R}$~: donc $R\in{\cal I}$.

\subsection{Fin de la d\'emonstration de la proposition \ref{proposition:stables}.}

Soit ${\cal J}$ un id\'eal premier de ${\cal N}$.
Si ${\cal J}\cap{\cal R}[q]=(0)$, alors ${\cal J}$ est principal et ne peut pas \^{e}tre $D$-stable
d'apr\`es la proposition
\ref{lemme:ideaux_principaux_hypergeo} (il ne peut contenir ni $q$ ni l'un des polyn\^{o}mes 
$y_0-y_1,y_0-y_2,y_1-y_2$). 

Supposons que ${\cal J}$ soit $D$-stable. Alors ${\cal Y}:={\cal J}\cap{\cal R}[q]\not=(0)$.
On a que ${\cal Y}$ est un id\'eal premier $D$-stable de ${\cal R}[q]$. 

Si ${\cal Y}\cap{\cal R}=(0)$ alors ${\cal Y}$ est principal, donc ${\cal Y}{\cal N}$
est un id\'eal premier principal $D$-stable de ${\cal N}$, et la proposition 
\ref{lemme:ideaux_principaux_hypergeo}
implique $q\in{\cal J}$. 

Sinon, ${\cal I}={\cal Y}\cap{\cal R}\not=(0)$ est un autre id\'eal 
premier $D$-stable, cette fois dans ${\cal R}$. Si ${\cal I}$ est principal alors ${\cal I}{\cal N}$
est principal et $D$-stable, et une autre application de la proposition \ref{lemme:ideaux_principaux_hypergeo} 
implique que l'un des polyn\^{o}mes $y_0-y_1,y_0-y_2,y_1-y_2$ appartient \`a ${\cal I}$.

Si ${\cal I}$ n'est pas principal alors, d'apr\`es  le
lemme \ref{lemme:nuovo}, $\tilde{{\cal I}}$ est un id\'eal premier non nul isobare et $D$-stable
de ${\cal R}$, et il contient un des polyn\^{o}mes $y_0-y_1,y_0-y_2,y_1-y_2$ d'apr\`es le lemme 
\ref{lemme:une_variable_hypergeo}.

\medskip

Il nous reste \`a justifier la
non nullit\'e de $\eta$.
\begin{Lemme} On a $\eta\not=0$.\label{lemme:elementaire}\end{Lemme}
\noindent {\bf Esquisse de d\'emonstration.} On doit utiliser les hypoth\`eses faites sur les rationnels $\alpha,\beta,\gamma$.
Analysons les facteurs d\'efinissant le produit $\eta$.
On a~:
\[a+b=1-(\alpha-\beta)^2.\] ce nombre est nul si et seulement si $\alpha=\beta-1$ ou $\alpha=\beta+1$, mais
ceci est clairement impossible car $1>\beta>\alpha>0$.
On a~:
\[a+c=\gamma(2-\gamma).\] ce nombre est clairement non nul, car $0<\gamma<1$. 
Le troisi\`eme facteur est~:
\[b+c=1-(\alpha+\beta)^2+\gamma(2(\alpha+\beta)-\gamma).\] Ce nombre est nul si et seulement si
$\gamma=-1+\alpha+\beta$ ou $\gamma=1+\alpha+\beta$. Si $\gamma=1+\alpha+\beta$, alors $\gamma>1$,
et cette eventualit\'e est exclue.
Si $\gamma=\alpha+\beta-1$, alors $\gamma-\alpha-\beta=-1$, mais $\gamma-\alpha-\beta>0$ par hypoth\`ese. Donc
ce troisi\`eme facteur est non nul.

Etudions, mais de mani\`ere esquiss\'ee, le quatri\`eme facteur~:
\begin{eqnarray*}
ab+ac+bc & = &
2(\gamma^2(-1+\alpha+\beta-2\alpha\beta)+\\
& &\gamma(1+\alpha+\beta)(1-\alpha-\beta+2\alpha\beta)-2(\alpha\beta)^2).\end{eqnarray*} On v\'erifie
directement que la fonction \[f(\alpha,\beta,\gamma)=ab+ac+bc\] est positive si
$\alpha,\beta,\gamma$ sont soumis aux contraintes $0<\alpha<\beta<\gamma<1$. 

La v\'erification de cette propri\'et\'e passe par une \'etude \'el\'ementaire assez calculatoire,
que nous ne reporterons pas ici. Les points critiques (c'est-\`a-dire les
points  qui annulent toutes les d\'eriv\'ees partielles de $f$ d'ordre $1$) de $f$
sont  donn\'es par les conditions simultan\'ees~:
\begin{eqnarray*}
(2\alpha-\gamma)(2\beta^2+\gamma-2\beta\gamma)& = & 0\\
(2\beta-\gamma)(2\alpha^2+\gamma-2\alpha\gamma)& = & 0\\
(1 - \alpha - \beta + 2 \alpha \beta) (1 + \alpha + \beta- 2 \gamma) & = & 0.
\end{eqnarray*}
On v\'erifie que les points critiques ne satisfont pas les in\'egalit\'es $0<\alpha<\beta<\gamma<1$
(\footnote{Il y a un seul point critique satisfaisant $0\leq\alpha\leq\beta\leq\gamma\leq 1$ qui 
est le point $(1/2,1/2,1)$. En calculant la matrice hessienne, on voit facilement que $f(1/2,1/2,1)=3/4$
est un maximum
r\'elatif.}). On passe 
ensuite \`a l'\'etude de la restriction de $f$ au bord du domaine de $\RR^3$ d\'efini par $0<\alpha<\beta<\gamma<1$.
Les points critiques des fonctions restriction sont alors encore donn\'es essentiellement par les conditions ci-dessus.
Apr\`es avoir examin\'e tous les cas, on pourra concl\^{u}re que $ab+ac+bc>0$.

\section{D\'emonstration du th\'eo\-r\`eme \ref{theorem:estimate2}.\label{section:demo_th}}

Soit ${\cal O}$ un ouvert non vide de $\CC$, soient $F_1(t),\ldots,F_m(t)$ des fonctions alg\'ebriquement ind\'ependantes sur
$K=\CC(t)$, admettant en tout point $\xi\in{\cal O}$ un d\'eveloppement en s\'erie de Puiseux convergent~:
\begin{equation}
F_i(t)=\sum_{k=s}^\infty v_{i,k}(t-\xi)^{k/q},\label{eq:Puiseux_F}\end{equation} avec $s\in\ZZ$, et pour un certain entier
$q>0$ ne d\'ependant pas de
$\xi$.

Soit $P_0$ un polyn\^{o}me de $\CC[t]$ non nul, et supposons que
l'anneau ${\cal A}=\CC[F_1,\ldots,F_m]$ soit stable pour la d\'erivation $\delta =P_0(t)\displaystyle{\frac{d}{dt}}$.
Il existe alors des polyn\^{o}mes $P_1,\ldots,P_m\in\CC[X_1,\ldots,X_m]$ tels que~:
\begin{eqnarray*}
\delta F_i & = & P_i(F_1,\ldots,F_m).
\end{eqnarray*}
Soit $d$ le plus grand degr\'e total des $P_i$. Soit $X_0$ une nouvelle inconnue, notons 
${\cal B}=K[X_0,\ldots,X_m]$ posons~:
\[B_i(X_0,\ldots,X_m) = X_0^dP_i\left(\frac{X_1}{X_0},\ldots,\frac{X_m}{X_0}\right)\in{\cal
B},\quad i=1,\ldots,m,\] et $B_0=X_0^dP_0(t)\in K[X_0]$. D\'efinissons l'op\'erateur de
d\'erivation homog\`ene~:
\[{\cal D}=B_0\frac{\partial }{\partial t}+X_0\sum_{i=1}^m B_i\frac{\partial}{ \partial X_i}.\]
On a donc, pour tout $Q\in{\cal B}$, \[({\cal D}Q)(1,F_1,\ldots,F_m)=\delta(Q(1,F_1,\ldots,F_m)).\]
De plus, si $Q\in{\cal B}$ est un \'el\'ement homog\`ene, alors ${\cal D}Q$ est
aussi homog\`ene. Nous appelons l'anneau diff\'erentiel $({\cal B},{\cal D})$ l'{\em anneau diff\'erentiel associ\'e \`a} 
$({\cal A},\delta)$.

\medskip

\noindent {\bf D\'efinition}. L'anneau diff\'erentiel $({\cal B},{\cal D})$ satisfait la {\em propri\'et\'e
de Ramanujan} s'il existe un polyn\^{o}me homog\`ene non nul $l\in{\cal B}$ tel que 
pour tout id\'eal premier non nul ${\mathfrak P}$ de ${\cal B}$ qui est ${\cal D}$-stable,
on a $l\in{\mathfrak P}$.

\medskip

Fixons un \'el\'ement $\xi\in{\cal O}$ et soit ${\cal K}$ le compl\'et\'e d'une cl\^{o}ture alg\'ebrique du compl\'et\'e de $K$ \`a la place
d\'etermin\'ee par $\xi$, soit ${\mathfrak P}$ un id\'eal premier homog\`ene de ${\cal B}$ et 
${\cal V}=Z({\mathfrak P})$ la sous-vari\'et\'e de $\PP_m({\cal K})$ des z\'eros de ${\mathfrak P}$. 
Notons 
$|\cdot|$ la valeur absolue de ${\cal K}$~:
si $F\in{\cal K}$, alors~:
\begin{equation}|F|=\exp\{-\mbox{ord}_\xi(F)\}.\label{eq:norma}\end{equation}

Nous posons
\[\overline{\omega}=(1:F_1:\cdots:F_m)\in\PP_m({\cal K}).\] Notons
\[\mbox{Dist}(\overline{\omega},{\cal V})=\mbox{Dist}_{(1,\ldots,1)}(\overline{\omega},{\cal V})\] 
la {\em distance de $\overline{\omega}$ \`a ${\cal V}$} de Philippon (voir p. 88 de
\cite{Nesterenko:Introduction6}), associ\'ee \`a la valeur absolue de ${\cal K}$. Cette  distance s'\'etend par multiplicativit\'e aux
cycles alg\'ebriques \'equidimen\-sionnels de
$\PP_m(K)$ (\footnote{ Tous les cycles alg\'ebriques seront positifs et \'equidimensionnels~: nous n'avons donc
plus \`a pr\'eciser systematiquement cet attribut.})~; elle d\'epend de $\xi$, et nous utiliserons parfois
la notation $\mbox{Dist}_{\xi}(\overline{\omega},{\cal V})$ lorsque la d\'ependance en $\xi$ doit \^{e}tre 
mise en relief. 

Posons \[\|\overline{\omega}\|=\max\{1,|F_1|,\ldots,|F_m|\},\] et pour un polyn\^{o}me
homog\`ene $U\in K[X_0,X_1,\ldots,X_m]$, 
\begin{equation}
U=\sum_{\und{\lambda}}c_{\und{\lambda}}\und{X}^{\und{\lambda}},
\label{eq:U}\end{equation}
posons aussi~:
\[|U|=\max_{\und{\lambda}}
\{|c_{\und{\lambda}}|\}.\]

On peut calculer facilement la distance de $\overline{\omega}$ \`a ${\cal V}$ si ${\cal V}$ est une hypersurface.
Soit $U$ comme dans (\ref{eq:U}), homog\`ene, avec $c_{\und{\lambda}}\in\CC[t]$, et supposons que ${\cal V}={\cal Z}(U)$.
Les
formules p. 64 de
\cite{Jadot:Criteres} impliquent~:
\begin{eqnarray}
\log\mbox{Dist}(\overline{\omega},{\cal V})& = & \log|U(\overline{\omega})|-\log|U|-
\deg_{\und{X}}(U)\log\|\overline{\omega}\|\label{eq:jatotto}.
\end{eqnarray}
En d'autres termes,
\begin{eqnarray*}
\log\mbox{Dist}(\overline{\omega},{\cal V})& =
&
-\mbox{ord}_{\xi}(U(\overline{\omega}))+\min_{\und{\lambda}}
\{\mbox{ord}_{\xi}(c_{\und{\lambda}})\}+\\ & &\deg(U)
\min_{i=1,\ldots,m}\{0,\mbox{ord}_{\xi}(F_i)\}.
\end{eqnarray*}

On rappelle \'egalement que l'image de l'application Dist$(\overline{\omega},\cdot)$ est contenue
dans l'intervalle $[0,1]$ (\cite{Nesterenko:Introduction7}, proposition 4.1 p. 119).

De plus, comme $F_1,\ldots,F_m$ satisfont (\ref{eq:Puiseux_F}) dans ${\cal O}$, on a
que pour tout cycle alg\'ebrique non nul ${\cal Z}$ de $\PP_m(K)$ (diff\'erent de $\PP_m(K)$), 
le nombre \[-\log\mbox{Dist}(\overline{\omega},
{\cal Z})\] est un \'el\'ement positif de $\ZZ/q$~; ceci se voit directement dans (\ref{eq:jatotto}) si ${\cal Z}=Z(Q)$
avec $Q\in {\cal B}-\{0\}$. Par
convention, nous posons
$-\log\mbox{Dist}(\overline{\omega}, Z(0)):=+\infty$ pour tout $\xi\in{\cal O}$.

Par convention, le cycle nul $0$ est de dimension $-1$, et $-\log\mbox{Dist}(\overline{\omega},
0)=0$ pour tout $\xi$.

Nous devons introduire une autre propri\'et\'e d'anneaux diff\'erentiels,
que nous appellons \og propri\'et\'e $D$\fg . 

\medskip

\noindent {\bf D\'efinition.} On dit que l'anneau diff\'erentiel $({\cal B},{\cal D} )$ satisfait la {\em propri\'et\'e
$D$ en $\xi$ d'exposant $\leq \rho(\xi)$}, si la constante $\rho(\xi)\geq 0$ ne d\'epend que de
$\xi$, et est telle que 
tout id\'eal premier homog\`ene ${\mathfrak P}$ de ${\cal B}$ contenant un id\'eal premier 
non nul $D$-stable, on ait~: 
\[-\log\mbox{Dist}(\overline{\omega},Z({\mathfrak P}))\leq\rho(\xi)
\deg({\mathfrak P}).\]
S'il existe une constante $\rho\geq 0$ telle
que pour tout $\xi\in{\cal O}$ l'anneau $({\cal B},{\cal D} )$ satisfasse la propri\'et\'e 
$D$ en $\xi$ d'exposant $\leq\rho$, alors nous dirons que $({\cal B},{\cal D} )$ satisfait
la {\em propri\'et\'e $D$ uniforme sur ${\cal O}$}.

\subsection{La propri\'et\'e de Ramanujan et la propri\'et\'e $D$.}

\begin{Lemme} Supposons que l'anneau diff\'erentiel $({\cal B},{\cal D})$ satisfasse la propri\'et\'e de 
Ramanujan. Alors il satisfait la propri\'et\'e $D$ en tout point $\xi\in{\cal O}$.
Soit $l$ le polyn\^{o}me fournit par la propri\'et\'e de Ramanujan.
Si~:
\[\rho=\sup_{\xi\in{\cal O}}\{-\log\mbox{Dist}_\xi(\overline{\omega},Z(l))\}\in\frac{\NN}{q},\]
alors $({\cal B},{\cal D})$ satisfait la propri\'et\'e $D$ uniforme sur ${\cal O}$ (d'exposant $\leq \rho$).
\label{lemme:distance}\end{Lemme}
\noindent {\bf D\'emonstration}. C'est l'analogue du lemme 3.4 p. 155 de \cite{Nesterenko:Introduction10}.
Fixons $\xi\in{\cal O}$, posons~:
\[\rho(\xi)=-\log\mbox{Dist}_\xi(\overline{\omega},Z(l)).\]
Sous les hypoth\`eses du lemme, $({\cal B},{\cal D})$ satisfait 
la propri\'et\'e $D$ d'exposant $\leq \rho(\xi)$~:
voici la preuve.

Soit ${\mathfrak P}$ un id\'eal premier homog\`ene non nul de ${\cal B}$, 
supposons que~:
\[-\log\mbox{Dist}(\overline{\omega},Z({\mathfrak P}))>\rho(\xi)
\deg({\mathfrak P}).\]
Soit $p$ un \'el\'ement homog\`ene non
nul de ${\mathfrak P}$. L'inegalit\'e suivante se d\'emontre 
en suivant les arguments du paragraphe 5 de \cite{Nesterenko:Introduction6},
ou en combinant les propositions 3.24 et 3.25 de \cite{Bosser:Independance}~:
\begin{eqnarray}
-\log\mbox{Dist}(\overline{\omega},Z(p)) & \geq-\displaystyle{\frac{1}{\deg({\mathfrak P})}}
\log\mbox{Dist}(\overline{\omega},Z({\mathfrak P}))
\label{eq:inegalite1}.
\end{eqnarray}
L'in\'egalit\'e (\ref{eq:inegalite1}) implique que $-\log\mbox{Dist}(\overline{\omega},
Z(p))>\rho(\xi)$, donc $p\not=l$. Ainsi $l\not\in{\mathfrak P}$, et d'apr\`es la propri\'et\'e de Ramanujan,
${\mathfrak P}$ est ne contient pas d'id\'eal premier non nul ${\cal D}$-stable.

Faisons varier $\xi\in{\cal O}$. Si le supremum $\rho$ est fini, il est clair
que $({\cal B},{\cal D})$ satisfait 
la propri\'et\'e $D$ uniforme d'exposant $\leq \rho$.

\medskip

Nous d\'emontrons maintenant~:
\begin{Proposition} Si ${\cal A}=\CC[t,e^t,Y_0,Y_1,Y_2]$ est muni de la d\'erivation $\delta$, alors
l'anneau diff\'erentiel homog\`ene associ\'e $({\cal B},{\cal D})$ satisfait la 
propri\'et\'e de Ramanujan, et la propri\'et\'e $D$ uniforme.
\label{proposition:prop_ramanujan}\end{Proposition}
\noindent {\bf D\'emonstration}.
Nous avons
${\cal B}=K[X_0,X_1,\ldots,X_4]$, $K=\CC(t)$ et 
\[{\cal D}=wX_0^2\frac{\partial}{\partial t}+wX_0^2X_1\frac{\partial}{\partial
X_1}+X_0\sum_{i=2}^4(X_i^2-U)\frac{\partial}{\partial X_i},\] 
avec 
\[U=\frac{1}{4}(a(X_2-X_3)^2+b(X_2-X_4)^2+c(X_3-X_4)^2).\]
Posons~:
\[{\cal A}^\sharp=K[e^t,Y_0(t),Y_1(t),Y_2(t)],\quad\overline{\omega}=(1:e^t:Y_0(t):
Y_1(t):Y_2(t))\in\PP_4({\cal K}).\]
on a un morphisme $\psi:{\cal B}\rightarrow{\cal A}^\sharp$ qui associe \`a
tout polyn\^{o}me homog\`ene $P\in{\cal B}$ l'\'el\'ement $\psi(P)=P(\overline{\omega})\in{\cal A}^\sharp$.
Clairement~: \[\psi({\cal D}P)=\delta(\psi (P)),\] pour tout polyn\^{o}me homog\`ene $P\in{\cal B}$, d'apr\`es
le lemme \ref{lemme:relations_explicites}.

D\'emontrons que si
${\mathfrak P}$ est un id\'eal premier homog\`ene non nul de ${\cal B}$ qui est ${\cal D}$-stable, 
alors $l\in{\mathfrak P}$ avec~:
\[l=X_0X_1(X_3-X_2)(X_4-X_2)(X_4-X_3).\]
Notons ${\cal P}^\sharp$ l'id\'eal de ${\cal A}^\sharp$ engendr\'e par les \'el\'ements de
$\psi({\mathfrak P})$~; c'est un id\'eal premier $\delta$-stable non nul.
L'id\'eal ${\cal P}={\cal P}^\sharp\cap{\cal A}$ de ${\cal A}$ est un id\'eal premier non nul $\delta$-stable
de ${\cal A}$.

Si ${\cal P}={\cal A}$, alors $X_0\in{\mathfrak P}$ et $l\in{\mathfrak P}$.
Si ${\cal P}\not={\cal A}$, au moins une des
fonctions $e^t,Y_2(t)-Y_0(t),Y_1(t)-Y_0(t),Y_2(t)-Y_1(t)$ appartient \`a ${\cal P}$, d'apr\`es la proposition
\ref{proposition:stables}. On en d\'eduit qu'au moins 
un des polyn\^{o}mes homog\`enes $X_1,X_3-X_2,X_4-X_2,X_4-X_3$ appartient \`a ${\mathfrak P}$~; finalement,
$l\in{\mathfrak P}$. Donc l'anneau diff\'erentiel $({\cal B},{\cal D})$ satisfait la propri\'et\'e $D$ d'apr\`es
le lemme \ref{lemme:distance}.

Le lemme
\ref{lemme:relations_explicites} implique~:
\begin{eqnarray*}
\delta\kappa & = & \delta(q(Y_0-Y_1)(Y_0-Y_2)(Y_1-Y_2))\\
& = & (q+2(Y_0+Y_1+Y_2))\kappa.
\end{eqnarray*}
De plus, la fonction $\kappa=l(\overline{\omega})$ ne s'annule pas dans $B^*$, car $(\kappa)$ est 
un id\'eal principal $\delta$-stable de $({\cal A},\delta)$, et $Y_0,Y_1,Y_2$ sont holomorphes sur $B^*$. 

Le nombre $\rho$ du lemme \ref{lemme:distance} est bien d\'efini (le calcul du supremum
se r\'eduit \`a un calcul de maximum), et d'apr\`es ce m\^{e}me lemme,
$({\cal A},\delta)$ satisfait la propri\'et\'e $D$ uniforme (On peut m\^{e}me montrer que l'exposant est $\leq 1$).

\subsection{Une estimation de multiplicit\'e provenant de la propri\'et\'e $D$.} 
\begin{Proposition} 
Supposons que l'anneau $({\cal B},{\cal D})$
satisfasse la propri\'et\'e $D$ d'exposant $\leq \rho(\xi)$ en $\xi$.
Il existe une constante $c_3>0$, d\'ependant uniquement de $\rho(\xi)$, avec la propri\'et\'e suivante. 
Soit $P$ un polyn\^{o}me non nul de $\CC[Z,X_1,\ldots,X_m]$ et posons~:
\[N_0=\max\{1,\deg_{Z}P\},\quad N_1=\max\{1,\deg_{X_1,\ldots,X_m}P\}.\]
Alors la fonction $F(z)=P(z,F_1(z)\ldots,F_m(z))$ satisfait \[\mbox{ord}_{\xi}F(z)\leq c_3N_0N_1^m.\]
\label{proposition:general}\end{Proposition}
La d\'emonstration de cette proposition est essentiellement la m\^{e}me que 
celle du th\'eor\`eme 2.2 p. 154 de \cite{Nesterenko:Introduction10}.
Nous suivons essentiellement l'approche de \cite{Bosser:Independance}, qui consiste
\`a utiliser les formes r\'esultantes (formes de Chow), plut\^{o}t que les formes \'eliminantes (comme
dans \cite{Nesterenko:Introduction10})~: cette approche est plus 
flexible lorsque l'on travaille avec des cycles alg\'ebriques dans des espaces projectifs.

\medskip

On note $h:K^\times\rightarrow\RR_{\geq 0}$ la {\em hauteur logarithmique de Weil}, que l'on \'etend de la mani\`ere usuelle
\`a $\PP_s(K)$, puis aux sous-vari\'et\'es et aux cycles de $\PP_s(K)$ via les formes de Chow (cf. paragraphe 3
de \cite{Nesterenko:Introduction6}). On note $\deg({\cal Z})$ le {\em degr\'e} d'un cycle d\'efini sur $K$
(paragraphe 2 de loc. cit.). 

Nous appellerons {\em id\'eal ${\cal D}$-instable} tout id\'eal premier homog\`ene non nul ${\mathfrak P}$ ne contenant aucun id\'eal
premier non nul ${\cal D}$-stable.

La proposition suivante est essentiellement la proposition 4.11 de 
\cite{Bosser:Independance}, et peut \^{e}tre aussi d\'emontr\'ee en utilisant les arguments
de \cite{Nesterenko:Introduction10} pp. 154-159. 

\begin{Proposition}
Il existe une constante $c_4>0$, ne d\'ependant que de $({\cal B},{\cal D})$, 
mais pas de $\xi$, satisfaisant la propri\'et\'e suivante.
Soit ${\mathfrak P}$ un id\'eal premier homog\`ene non nul ${\cal D}$-instable de ${\cal B}$, soit $p\in{\mathfrak
P}\cap\CC[t][X_0,\ldots,X_m]$ un polyn\^{o}me homog\`ene non nul tel que la quantit\'e
\[\nu({\mathfrak P})=\max\{1,h({\mathfrak P})\}\deg_{\und{X}}(p)+\deg({\mathfrak P})\deg_t(p)\] soit minimale. Alors il existe
un entier
$N$ avec $1\leq N\leq c_4$, tel que pour tout $n<N$ on ait ${\cal D}^np\in{\mathfrak P}$, et tel que
${\cal D}^Np\not\in{\mathfrak P}$.
\label{proposition:411}\end{Proposition}

Les arguments p. 158 de
\cite{Nesterenko:Introduction10} permettent de montrer qu'un id\'eal premier non nul ${\mathfrak P}$ de ${\cal B}$ est 
${\cal D}$-instable si et seulement si pour tout id\'eal non nul ${\mathfrak I}\subset{\mathfrak P}$ on a
${\cal D}{\mathfrak I}\not\subset{\mathfrak I}$. 
En particulier, pour tout $p\in{\mathfrak P}-\{0\}$ il existe
un entier $n$ tel que ${\cal D}^np\not\in{\mathfrak P}$~: la proposition \ref{proposition:411}
fournit un choix privil\'egi\'e pour $p$, de telle mani\`ere que son degr\'e soit minimal. 

L'in\'egalit\'e (34), ou l'in\'egalit\'e (35) p. 144 de \cite{Nesterenko:Introduction9}, ou encore le
lemme 3.1 p. 154 de \cite{Nesterenko:Introduction10} (ce sont essentiellement des majorations explicites de fonctions de
Hilbert bi-homog\`enes), impliquent que $p$ satisfait~:
\begin{eqnarray}
\deg(p) & \leq & c_5\deg({\mathfrak P})^{1/(m-g)}\nonumber\\
h(p) & = & \deg_t(p)\label{eq:egalite_irred}\\
& \leq & c_5\max\{1,h({\mathfrak P})\}\deg({\mathfrak P})^{-(m-g-1)/(m-g)}\nonumber,
\end{eqnarray}
o\`u $g=\dim({\mathfrak P})$ et $c_5$ est une constante positive~: on a l'\'egalit\'e (\ref{eq:egalite_irred}) car $p$ est
irr\'eductible dans $\CC[t,X_0,\ldots,X_m]$, d'apr\`es la minimalit\'e 
de $\nu({\mathfrak P})$. La d\'efinition de la d\'erivation ${\cal D}$ implique
que $\deg({\cal D}^Np)\leq\deg(p)+Nd\leq (d+1)c_4\deg(p)$, et $h({\cal D}^Np)\leq\deg_t({\cal D}^Np)\leq
h(p)+Nd\leq (d+1)c_4 h(p)$.
Donc~:
\begin{eqnarray}
\deg({\cal D}^Np) & \leq & c_6\deg({\mathfrak P})^{1/(m-g)}\label{eq:degre}\\
h({\cal D}^Np) & \leq & 
c_6\max\{1,h({\mathfrak P})\}\deg({\mathfrak P})^{-(m-g-1)/(m-g)}\label{eq:hauteur},
\end{eqnarray}
avec $c_6=(d+1)c_4$.
Nous utilisons toutes ces propri\'et\'es dans la proposition
suivante.

\begin{Proposition}
Supposons que l'anneau diff\'erentiel $({\cal B},{\cal D})$ satisfasse la propri\'et\'e $D$ en $\xi$ 
d'exposant $\leq \rho(\xi)$.

Il existe deux constantes $c_7,c_{8}>0$, ne d\'ependant que de 
$\rho(\xi)$, telles que la propri\'et\'e suivante soit satisfaite.
Soit ${\mathfrak P}$ un id\'eal premier homog\`ene non nul de ${\cal B}$, soit $g=\dim({\mathfrak P})$, supposons
que 
\begin{equation}
\log\mbox{{\em Dist}}(\overline{\omega},Z({\mathfrak P}))\leq -c_7(\max\{1,h({\mathfrak P})\}\deg({\mathfrak P})^{(g+1)/(m-g)}+\deg({\mathfrak
P})^{m/(m-g)}).
\label{eq:hypothese_contradictoire}\end{equation}
 Alors il existe un polyn\^{o}me homog\`ene $\phi\in K[X_0,\ldots,X_m]\cap \CC[t,X_0,\ldots,X_m]$ n'appartenant pas \`a ${\mathfrak P}$,
satisfaisant
\begin{eqnarray}
\deg_{\und{X}}(\phi) & \leq & c_8\deg({\mathfrak P})^{1/(m-g)}\label{eq:degre1}\\
h(\phi) & \leq & c_8\max\{1,h({\mathfrak P})\}\deg({\mathfrak P})^{-(m-g-1)/(m-g)}\label{eq:hauteur1},
\end{eqnarray}
et tel que pour tout $\overline{v}\in \PP_m({\cal K})$, $\overline{v}\in Z({\mathfrak P})$ on ait~:
\begin{equation}\mbox{{\em Dist}}(\overline{\omega},Z(\phi))\leq\mbox{{\em Dist}}(\overline{\omega},\{\overline{v}\}).
\label{eq:distance_proposition}\end{equation}
\label{proposition:existence_de_polynome}\end{Proposition}
\noindent {\bf D\'emonstration.} Ceci corr\'espond \`a la proposition 3.6 p. 156 de \cite{Nesterenko:Introduction10}, ou
\`a la proposition 4.18 de \cite{Bosser:Independance}, mais nous 
donnons les d\'etails de la d\'emonstration, car nous avons suppos\'e que les
fonctions $F_1,\ldots,F_m$ satisfont (\ref{eq:Puiseux_F}), alors que dans ces
derni\`eres r\'ef\'erences les fonctions sont toutes holomorphes en $\xi$.

Si l'in\'egalit\'e
(\ref{eq:hypothese_contradictoire}) est satisfaite pour une certaine constante $c_7$, 
alors on a aussi les deux in\'egalit\'es~:
\begin{eqnarray}
\log\mbox{Dist}(\overline{\omega},Z({\mathfrak P})) & \leq & -c_7(\max\{1,h({\mathfrak P})\}\deg({\mathfrak P})^{1/(m-g)}+\deg({\mathfrak
P}))\\
& \leq & -c_7\deg({\mathfrak P})\nonumber.
\end{eqnarray}
Si $c_7\geq \rho(\xi)$, alors l'id\'eal ${\mathfrak P}$ est ${\cal D}$-instable (propri\'et\'e $D$), et la proposition 
\ref{proposition:411} peut \^{e}tre appliqu\'ee. Soient donc $p,N$ comme dans cette proposition.
On a que $\tilde{p}={\cal D}^{N-1}p\in{\mathfrak P}$ d'apr\`es la minimalit\'e de $N$, et
$-\log$ Dist$(\overline{\omega},Z(\tilde{p}))\geq c_7$, d'apr\`es l'in\'egalit\'e (\ref{eq:inegalite1}).
Posons aussi $\tilde{q}={\cal D}\tilde{p}={\cal D}^Np$.
On a
$\tilde{q}\not\in{\mathfrak P}$, $\tilde{q}\in\CC[t][X_0,\ldots,X_m]$, et~:
\begin{eqnarray*}
\tilde{q}(\overline{\omega}) & = & ({\cal D}\tilde{p})(\overline{\omega})\\
& = & \delta (\tilde{p}(\overline{\omega}))\\
& = & P_0(t)\frac{d}{dt}\tilde{p}(\overline{\omega}),
\end{eqnarray*}
d'apr\`es la d\'efinition de ${\cal D}$. Donc, puisque $P_0$ est un polyn\^{o}me~:
\[\mbox{ord}_{\xi}(\tilde{q}(\overline{\omega}))\geq\mbox{ord}_{\xi}(\tilde{p}(\overline{\omega}))-1.\]
En posant $U=\tilde{q}$ dans l'in\'egalit\'e (\ref{eq:jatotto}) on obtient~:
\begin{eqnarray*}
-\log\mbox{Dist}(\overline{\omega},Z(\tilde{q})) & \geq & -\log|\tilde{p}(\overline{\omega})|
+\log|\tilde{p}|+\deg(\tilde{p})\log\|\overline{\omega}\|-\\
& & \log|\tilde{p}|-\deg(\tilde{p})\log\|\overline{\omega}\|+\\
& & \log|\tilde{q}|+\deg(\tilde{q})\log\|\overline{\omega}\|-1\\
& \geq & -\log\mbox{Dist}(\overline{\omega},Z(\tilde{p}))+\log|\tilde{q}|-\log|\tilde{p}|+\\
& & (\deg(\tilde{q})-\deg(\tilde{p}))\log\|\overline{\omega}\|-1\\
& \geq & -\log\mbox{Dist}(\overline{\omega},Z(\tilde{p}))+\log|\tilde{q}|+\\
& & d\log\|\overline{\omega}\|-1\\
& \geq & -\log\mbox{Dist}(\overline{\omega},Z(\tilde{p}))+\log|\tilde{q}|-1,
\end{eqnarray*}
(noter que $\log\|\overline{\omega}\|$ est un nombre positif ou nul). 

On montre facilement que \[\log|\tilde{q}|-1\geq (1/2)\log\mbox{Dist}(\overline{\omega},Z(\tilde{p}))\]
(on utilise les 
m\^{e}mes arguments p. 159 de \cite{Nesterenko:Introduction10}), et on parvient \`a l'in\'egalit\'e
\[-\log\mbox{Dist}(\overline{\omega},Z(\tilde{q}))\geq -(1/2)\log\mbox{Dist}(\overline{\omega},Z(\tilde{p})).\]

Soit $\overline{v}\in Z({\mathfrak P})\subset\PP_m({\cal K})$. 
D'apr\`es le corollaire 4.9 p. 42 de \cite{Nesterenko:Introduction3},ou la proposition
3.25 de \cite{Bosser:Independance}, on a que \[-\log\mbox{Dist}(\overline{\omega},Z(\tilde{p}))\geq-\log\mbox{Dist}
(\overline{\omega},\{\overline{v}\}).\]
Posons $\phi=\tilde{q}^2$. On a $\phi\not\in{\mathfrak P}$ et~:
\begin{eqnarray*}
-\log\mbox{Dist}(\overline{\omega},Z(\phi)) & = & -2\log\mbox{Dist}(\overline{\omega},Z(\tilde{q}))\\
& \geq & -\log\mbox{Dist}(\overline{\omega},Z(\tilde{p}))\\
& \geq & -\log\mbox{Dist}(\overline{\omega},\{\overline{v}\}).
\end{eqnarray*}
D'autre part, d'apr\`es les in\'egalit\'es (\ref{eq:degre}) et (\ref{eq:hauteur}),
on a $\deg(\phi)\leq 2c_6\deg({\cal D}^Np)$ et $h(\phi)\leq 2c_6h({\cal D}^Np)$. En posant $c_8=2c_6$, on
voit que les in\'egalit\'es (\ref{eq:degre1}) et (\ref{eq:hauteur1}) sont satisfaites, et toutes les
constantes ne d\'ependent que de $\rho(\xi)$, ce qui 
compl\`ete la d\'emonstration de la proposition \ref{proposition:existence_de_polynome}.

\medskip

Nous d\'emontrons la proposition \ref{proposition:general}, en d\'emontrant le r\'esultat encore plus g\'en\'eral
suivant (cf. th\'eor\`eme 2.2 p. 154 de \cite{Nesterenko:Introduction10}).
\begin{Proposition}
Supposons que l'anneau diff\'erentiel $({\cal B},{\cal D})$ satisfasse la propri\'et\'e
$D$ en $\xi$ d'exposant $\leq\rho(\xi)$.
Alors il existe une constante $c_9>0$, ne d\'ependant que de $\rho(\xi)$, telle que pour tout cycle alg\'ebrique 
${\cal Z}\subset\PP_m({\cal K})$ de dimension $g<m$, 
l'in\'egalit\'e suivante soit satisfaite~:
\begin{equation}
\log\mbox{{\em Dist}}(\overline{\omega},{\cal Z})\geq -c_9(\max\{1,h({\cal Z})\}\deg({\cal Z})^{(g+1)/(m-g)}+
\deg({\cal Z})^{m/(m-g)}).\label{eq:step0}\end{equation}
\label{proposition:encore_plus_general}\end{Proposition}
\noindent {\bf D\'emonstration}. 
On fait une d\'emonstration par r\'ecurrence sur $g$. Pour $g=-1$, la proposition est trivialement vraie
car ${\cal Z}$ est le cycle nul.
En suivant des arguments similaires \`a ceux pp. 159-160 de \cite{Nesterenko:Introduction10}, ou 
en suivant \cite{Bosser:Independance} (on utilise les propri\'et\'es d'additivit\'e de $h(\cdot)$ et $\deg(\cdot)$
sur les cycles), on
suppose que la proposition \ref{proposition:encore_plus_general} soit vraie pour tout cycle 
${\cal Z}$ de dimension $g\geq -1$, et nous la d\'emontrons uniquement pour une vari\'et\'e
${\cal V}\subset\PP_m({\cal K})$ de dimension $g+1$. Supposons par l'absurde que la proposition 
\ref{proposition:encore_plus_general} soit fausse pour les vari\'et\'es de dimension $g+1$ de $\PP_m({\cal K})$.
Alors, pour toute constante $c_9>0$ assez grande, il existe un id\'eal premier homog\`ene non nul ${\mathfrak P}$ de dimension
$g+1$ de ${\cal B}$
tel que \begin{small}
\[\log\mbox{Dist}(\overline{\omega},Z({\mathfrak P}))< -c_9(\max\{1,h({\mathfrak P})\}\deg({\mathfrak P})^{(g+2)/(m-g-1)}+
\deg({\mathfrak P})^{m/(m-g-1)}).\]\end{small}
Nous pouvons choisir $c_9>c_7$ et travailler avec un id\'eal premier ${\mathfrak P}$ (associ\'e au choix de
la constante $c_9$) qui
satisfait les conditions de la proposition \ref{proposition:existence_de_polynome}. Soit $\phi$ le polyn\^{o}me
fournit par cette proposition, consid\'erons l'id\'eal ${\mathfrak Q}=({\mathfrak P},w)$~:
il est pur de dimension $g$. Soit ${\cal Z}$ le cycle associ\'e~: nous pouvons appliquer 
les th\'eor\`emes de B\'ezout g\'eom\'etrique et arithm\'etique (cf. paragraphe 4 de \cite{Nesterenko:Introduction6}
pour des versions archim\'ediennes).
\begin{eqnarray*}
\deg({\cal Z}) & \leq & \deg({\mathfrak P})\deg(\phi)\\
& \leq &  c_9\deg({\mathfrak P})^{1+1/(m-g-1)}\\
& \leq & c_9\deg({\mathfrak P})^{(m-g)/(m-g-1)},\\
h({\cal Z}) & \leq & h({\mathfrak P})\deg(\phi)+h(\phi)\deg({\mathfrak P})\\
& \leq & c_9(h({\mathfrak P})\deg({\mathfrak P})^{1/(m-g-1)}+\\
& & \max\{1,h({\mathfrak P})\}\deg({\mathfrak
P})^{1-(m-g-2)/(m-g-1)})\\ & \leq & 2c_9\max\{1,h({\mathfrak P})\}\deg({\mathfrak P})^{1/(m-g-1)}.
\end{eqnarray*}
Donc nous avons l'in\'egalit\'e~:
\begin{eqnarray}
\lefteqn{\max\{1,h({\cal Z})\}\deg({\cal Z})^{(g+1)/(m-g)}+\deg({\cal Z})^{m/(m-g)}\leq}\nonumber\\
& \leq & 2c_9\max\{1,h({\mathfrak P})\}\deg({\mathfrak P})^{1/(m-g-1)}\times\nonumber\\ & &
\left(c_9\deg({\mathfrak
P})^{(m-g)/(m-g-1)}\right)^{(g+1)/(m-g)}+\nonumber\\ & & \left(c_9\deg({\mathfrak
P})^{(m-g)/(m-g-1)}\right)^{m/(m-g)}\nonumber\\ & \leq & 2c_9^{1+r/(m-g)}\max\{1,h({\mathfrak P})\}\deg({\mathfrak
P})^{(g+2)/(m-g-1)}+\nonumber\\ & & c_9^{m/(m-g)}
\deg({\mathfrak P})^{m/(m-g-1)}\nonumber\\
& \leq & c_{12}(\max\{1,h({\mathfrak P})\}\deg({\mathfrak P})^{(g+2)/(m-g-1)}+\deg({\mathfrak P})^{m/(m-g-1)})\label{eq:step1},
\end{eqnarray}
pour une certaine constante $c_{10}$, d\'ependant uniquement de $c_9$ et $m$. En combinant avec
l'in\'egalit\'e (\ref{eq:step0}), l'in\'egalit\'e (\ref{eq:step1}) implique~:
\begin{equation}
\log\mbox{Dist}(\overline{\omega},{\cal Z})\geq-c_{11}(\max\{1,h({\mathfrak P})\}\deg({\mathfrak
P})^{(g+2)/(m-g-1)}+\deg({\mathfrak P})^{m/(m-g-1)})\label{eq:step2}
\end{equation}
pour une constante $c_{11}>0$ d\'ependant de $c_6,c_9$ et $m$.

D'autre part, la condition m\'etrique (\ref{eq:distance_proposition}) permet d'appliquer 
le premier th\'eor\`eme de B\'ezout m\'etrique (cf. \cite{Nesterenko:Introduction8}, paragraphe 7 
pour une version archi\-m\'edienne). On a~:
\begin{eqnarray}
\log\mbox{Dist}(\overline{\omega},{\cal Z})& \leq & \log\mbox{Dist}(\overline{\omega},Z({\mathfrak P}))+h(\phi)\deg({\mathfrak
P})+\deg(\phi)h({\mathfrak P})\nonumber\\
& \leq & \log\mbox{Dist}(\overline{\omega},Z({\mathfrak P}))+2c_9\max\{1,h({\mathfrak P})\}\deg({\mathfrak P})^{1/(m-g-1)}\nonumber\\
& \leq & -c_{12}(h({\mathfrak P})\deg({\mathfrak P})^{(g+2)/(m-g-1)}+
\deg({\mathfrak P})^{m/(m-g-1)})+\nonumber\\
& &2c_9\max\{1,h({\mathfrak P})\}\deg({\mathfrak P})^{1/(m-g-1)}\label{eq:step3}.
\end{eqnarray}
En combinant les in\'egalit\'es (\ref{eq:step2}) et (\ref{eq:step3}) nous trouvons que, pour tout $c_{12}>c_7$,
il existe un id\'eal premier homog\`ene ${\mathfrak P}$ de dimension $g+1$, tel que~:
\begin{eqnarray*}
\lefteqn{c_{12}(\max\{1,h({\mathfrak P})\}\deg({\mathfrak P})^{(g+2)/(m-g-1)}+
\deg({\mathfrak P})^{m/(m-g-1)})\leq}\\& \leq & c_{13}(\max\{1,h({\mathfrak P})\}\deg({\mathfrak P})^{(g+2)/(m-g-1)}+
\deg({\mathfrak P})^{m/(m-g-1)}),\end{eqnarray*} pour une constante $c_{13}>0$ qui elle, ne d\'epend pas
de $c_{12}$, d'o\`u une contradiction.

\medskip

Voici comment la proposition \ref{proposition:encore_plus_general} implique la proposition
\ref{proposition:general}. Soit $P\in\CC[Z,X_1,\ldots,X_m]$: alors $P\in K[X_1,\ldots,X_m]$. Soit $U\in{\cal B}$
l'homog\'eneis\'e de $P$ et posons ${\cal Z}=Z(U)\subset\PP_m(K)$. Nous pouvons supposer que $\deg_{\underline{X}}(P)\not=0$.
On a~:
\[h({\cal Z})\leq N_0,\quad\deg({\cal Z})\leq N_1,\quad g=m-1.\]
La proposition \ref{proposition:encore_plus_general} implique~:
\begin{eqnarray}
-\log\mbox{Dist}(\overline{\omega},{\cal Z}) & \leq & c_{12}(N_0N_1^m+N_1^m)\nonumber\\
& \leq & 2c_{12}N_0N_1^m\label{eq:joppa}
\end{eqnarray}
Comme $\log|U|\geq -\deg_Z(U)$, on a~:
\begin{eqnarray*}
\lefteqn{\mbox{ord}_\xi(U(\overline{\omega}))+\log|U|+\deg_{\underline{X}}(U)\log\|\overline{\omega}\|\geq}\\
& \geq & \mbox{ord}_\xi(U(\overline{\omega}))-N_0.
\end{eqnarray*}
En utilisant l'\'egalit\'e (\ref{eq:jatotto}) et en combinant cette derni\`ere estimation avec
l'in\'egalit\'e (\ref{eq:joppa}), on trouve~:
\begin{eqnarray*}
\mbox{ord}_\xi(F) & = & \mbox{ord}_\xi(U(\overline{\omega}))\\
& \leq & 2c_{12}N_0N_1^m+N_0\\
& \leq & c_{13}N_0N_1^m.
\end{eqnarray*}
Le lecteur peut v\'erifier que toutes ces constantes ne d\'ependent que de $\rho(\xi)$. Le th\'eor\`eme
\ref{theorem:estimate2} est une cons\'equence de la proposition \ref{proposition:general},
gr\^{a}ce \`a la proposition \ref{proposition:prop_ramanujan}.

\end{document}